\documentclass[a4paper,11pt]{article}

\usepackage{amssymb,amsmath,amsthm,amsfonts}
\usepackage{mathrsfs}
\usepackage{multirow}
\usepackage{graphicx}
\usepackage{authblk}
\usepackage{indentfirst}
\usepackage{multicol}
\usepackage{tabu}
\usepackage{url}
\usepackage{fancyhdr}
\usepackage[numbers,sort&compress]{natbib}
\usepackage[all]{xy}
\usepackage{mathtools}
\usepackage[english]{babel}
\usepackage{colortbl}
\usepackage{caption}
\usepackage{hyperref}
\usepackage{bookmark}
\usepackage{subcaption}
\usepackage{tikz}
\usepackage{float}
\usepackage{enumitem}
\usepackage{xcolor}
\usepackage{tikz-cd}
\usepackage{bm}
\usepackage{booktabs}
\usepackage{longtable}
\usepackage{extarrows}
\usepackage[section]{placeins}
\usepackage{array}
\usepackage{verbatim}
\usepackage{adjustbox}\usepackage[UTF8,fontset=fandol,scheme=plain]{ctex}
\usetikzlibrary{arrows.meta}

\newtheorem{theorem}{Theorem}[section]
\newtheorem{corollary}[theorem]{Corollary}

\newtheorem{lemma}[theorem]{Lemma}
\newtheorem{proposition}[theorem]{Proposition}

\theoremstyle{definition}
\newtheorem{definition}[theorem]{Definition}
\newtheorem{example}[theorem]{Example}
\newtheorem{remark}[theorem]{Remark}

\title{An explicit categorification of a Jones polynomial for tangles}

\author[1,2]{Jing-Wen Gao}
\author[1]{Jian Liu}
\author[1]{Kefeng Liu}
\author[3]{Li Shen}

\affil[1]{Mathematical Science Research Center, Chongqing University of Technology, Chongqing 400054, China}
\affil[2]{School of Arts and Sciences, Fujian Medical University, Fuzhou 350122, China}
\affil[3]{NSF-Simons National Institute for Theory and Mathematics in Biology, Chicago, IL 60637, USA.}

\makeatletter
    \renewcommand*{\@fnsymbol}[1]{\ensuremath{\ifcase#1\or \dagger\or *\or *\or
   \mathsection\or \else\@ctrerr\fi}}
\makeatother

\date{}

\begin{document}

    \maketitle

    \paragraph{Abstract}

    We define a Jones polynomial for oriented tangles and prove that it is categorified by their TQFT-based Khovanov homology. Extending the Kauffman bracket to general $n$-tangles and applying an augmentation that sends each circle-free flat tangle to $1$, we obtain an unnormalised Jones polynomial $\mathcal{J}(T)$ whose normalisation recovers the classical Jones polynomial for links. We establish a skein relation for $\mathcal{J}$ and prove the categorification identity
    \[
    \mathcal{J}(T)=q^{n}\chi(H(T)),
    \]
    where $\chi(H(T))$ denotes the graded Euler characteristic of the Khovanov homology of $T$. This identity extends the classical categorification formula for links. For a zero-framed single-arc tangle, we identify the suitably graded through-state subcomplex of its $(N-1)$-parallel cable with Khovanov's finite reduced colored complex. Its graded Euler characteristic is the reduced $N$-colored Jones polynomial of the closure. We also describe basic tangle operations, stacking formulas, and closure-volume applications of known colored Jones asymptotics.



\tableofcontents

\section{Introduction}

Knots and links are basic objects in low-dimensional topology. A tangle is a piece of a knot or link embedded in a $3$-ball with a specified set of boundary points. By decomposing links into tangles, local topological structures can be studied independently. Various invariants have been introduced to classify tangles. Unlike closed link invariants, tangle invariants must incorporate the boundary data. Consequently, the construction of tangle invariants relies on algebraic or categorical formalisms capable of encoding both the internal topology and the boundary constraints.

The Jones polynomial, introduced in \cite{Jones1985polynomial}, is a fundamental knot invariant. The two trefoils share the same Alexander polynomial~\cite{Alexander1928}, but their Jones polynomials differ, as Table~\ref{computation:at3} shows for $3_1$ and $3'_1$. The invariant is related to statistical mechanics, quantum groups, and quantum field theory, and it can be computed combinatorially via the Kauffman bracket \cite{Kauffman1987}. Khovanov later introduced a categorification of the Jones polynomial, lifting it to a bigraded homology theory \cite{khovanov2000categorification}.

The extension of the Jones polynomial to tangles has been developed along several lines. In \cite{turaev1990jones}, an isotopy invariant of tangles was introduced, generalizing the classical Jones polynomials. Furthermore, Khovanov \cite{khovanov2002functor} categorified the Jones polynomial for tangles by defining a functor-valued invariant that assigns a complex of bimodules to a tangle diagram. Recently, within the framework of representation theoretic categorified tangle invariants, \cite{stroppel2023categorification} assigns to each tangle a complex of Soergel bimodules whose induced morphism on Grothendieck groups recovers the Reshetikhin-Turaev invariant, yielding the Jones polynomial as the Euler characteristic upon closure.

In \cite{liu2026khovanov,liu2024persistent}, the authors provided a TQFT‑based construction that yields explicit computations of Khovanov homology for tangles. It is natural to ask whether the Khovanov homology obtained from this construction gives the categorification of the Jones polynomial for tangles. Addressing this question offers a more elementary, direct, and distinct perspective on the categorification of the Jones polynomial for tangles.

In this paper, we define a Jones polynomial for tangles and prove that it is categorified by the Khovanov homology constructed in \cite{liu2026khovanov,liu2024persistent}, thereby answering the question raised above affirmatively. We first define the Kauffman bracket for  $n$-tangles. Recall that Khovanov \cite{khovanov2002functor} defined a bracket for oriented $(m,n)$-tangles; we modify the normalisation so that, for links, our bracket coincides with the classical one \cite{Kauffman1987,khovanov2000categorification}. The resulting bracket takes values in the free $\mathbb{Z}[q,q^{-1}]$-module generated by the isotopy classes of circle-free flat $n$-tangles. The augmentation sends each such basis element to $1$ and gives the \emph{unnormalised Jones polynomial} $\mathcal{J}(T)$. For a nonempty link, the normalisation $\widetilde{\mathcal{J}}(T)=\mathcal{J}(T)/(q^{-1}+q)$ recovers the classical Jones polynomial under the substitution $t^{1/2}=-q$. We further show a skein relation for \(\mathcal{J}(T)\). Theorem~\ref{Main} generalizes Khovanov's classical formula \cite{khovanov2000categorification} from links to tangles. Specifically, for an oriented $n$-tangle $T$ with Khovanov homology $H(T)$, we prove
\[
\mathcal{J}(T)=q^{n}\,\chi(H(T)),
\]
which recovers the classical statement in the case $n=0$. We verify this identity by computing both invariants for tangles with at most three crossings. Finally, we examine the behaviour of $\mathcal{J}$ under three basic operations. It is multiplicative under disjoint union, acquires a factor of $(q^{-1}+q)$ for each loop under the gluing of $1$-tangles, and satisfies $\mathcal{J}(T^{*})(q)=\mathcal{J}(T)(q^{-1})$ under mirror reflection.

Section~\ref{sec:colored-application} connects the local tangle construction with Khovanov's finite reduced colored theory \cite{khovanov2002functor,khovanov2005colored}. For a zero-framed single-arc $(1,1)$-tangle $T$ with no closed components, Theorem~\ref{thm:app-comparison} identifies the suitably graded through-state subcomplex of its $(N-1)$-parallel cable with Khovanov's finite reduced colored complex. This comparison at the chain level yields diagram independence and gives
\[
\chi(H_N(T))=J_N(\widehat T;q^2),
\]
where $\widehat T$ is the standard closure. Theorem~\ref{thm:app-stacking} gives a chain-level tensor-product formula under stacking, implying multiplicativity of the Poincar\'e polynomials and graded Euler characteristics of the reduced colored homology. Combining this with known asymptotics of the colored Jones polynomial, we show that at $q=e^{\pi i/N}$ the logarithmic growth of these Euler characteristics  recovers the closure volume for single-arc representatives of the figure-eight and torus knots and their finite stacking products. Here the closure volume denotes the sum of the volumes of the hyperbolic pieces of the closure exterior.

This work is organized as follows. In the next section, we review tangles and their associated Khovanov homology. In Section \ref{section3}, we extend the Kauffman bracket to tangles, define the corresponding Jones polynomial, and establish its relation to the graded Euler characteristic of Khovanov homology for tangles. Section~\ref{section5} examines the behavior of the Jones polynomial under several fundamental topological operations. Section~\ref{sec:colored-application} gives the reduced colored extension, its comparison and stacking properties, and the closure-volume applications.

\section{Preliminaries}\label{sec:preliminaries}

In this section, we recall some basic notions and results concerning tangles. For a general introduction to the theory of tangles, we refer the reader to \cite{le1995representation}. The classical construction of Khovanov homology for tangles and related results can be found in \cite{bar2005khovanov,khovanov2002functor,liu2024persistent}.

\subsection{Tangles and states}

An $n$-tangle is a proper embedding of $n$ arcs together with a finite collection of circles in the $3$-ball $B^3$, with the endpoints of the arcs lying on the boundary $\partial B^3$. A \emph{flat $n$-tangle} \cite{khovanov2002functor} is an $n$-tangle that can be embedded in the $2$-ball $B^2$ so that the images of all arcs and circles are pairwise disjoint, with the endpoints of the arcs lying on the boundary $\partial B^2$.
The collection of circles is allowed to be empty. In particular, the empty tangle is regarded as a flat $0$-tangle.

Let $T$ be a tangle and $D$ be its diagram. A crossing of $D$ can be resolved in one of the following two ways.

\begin{figure}[htbp]
	\centering 
	\includegraphics[width=0.45\textwidth]{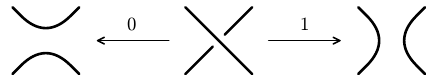}
	\caption{The \(0\)-resolution (left) and \(1\)-resolution (right) of the central crossing.}
	\label{resolution2}
\end{figure}

A \emph{state} $D_s$ of a diagram $D$ is a diagram in which each crossing of $D$ has been resolved in one of the two ways shown in Figure \ref{resolution2}. Thus, $D_s$ is a flat tangle. Clearly, a diagram with $c$ crossings has $2^c$ states. Label the crossings of $D$ by $1, 2, \ldots, c$. Each state is then determined by a sequence of $c$ entries, each of which is either $0$ or $1$. Thus, the resulting $2^c$ states are naturally indexed by the elements of $\{0,1\}^c$.

Let $I$ denote the set of crossings of $D$. To each subset $J \subseteq I$, we associate the state $D(J)$ obtained by taking the $1$-resolution at every crossing in $J$ and the $0$-resolution at every crossing in $I \setminus J$. We regard the subsets of $I$ as the vertices of a directed cube. Whenever $J \subset J'$ and $|J' \setminus J| = 1$, we connect $J$ to $J'$ by a directed edge. In this way, the collection of all states naturally forms a directed $c$-dimensional cube.

\subsection{TQFT construction of tangles}

 Throughout the paper, unless otherwise specified, all vector spaces are understood to be over the field of rational numbers $\mathbb{Q}$. Let
\[
V=\mathbb{Q}\{v_+,v_-\},\qquad
W=\mathbb{Q}\{w\},
\]
where
\[
\deg(v_+)=1, \quad \deg(v_-)=\deg(w)=-1.
\]
We define the following linear maps, which endow $V$ with a Frobenius algebra structure and $W$ with a Frobenius module structure:
\[
\begin{aligned}
m_1:\quad & V\otimes_{\mathbb{Q}}V \longrightarrow V, & \quad & \left\{ \begin{array}{ll} v_+\otimes v_+ \longmapsto v_+, & v_-\otimes v_+ \longmapsto v_-, \\ v_+\otimes v_- \longmapsto v_-, & v_-\otimes v_- \longmapsto 0 \end{array} \right. \\[2mm]
m_2:\quad & W\otimes_{\mathbb{Q}}V \longrightarrow W, & \quad & \left\{ \begin{array}{l} w\otimes v_+ \longmapsto w, \\ w\otimes v_- \longmapsto 0 \end{array} \right. \\[2mm]
\Delta_1:\quad & V \longrightarrow V\otimes_{\mathbb{Q}}V, & \quad & \left\{ \begin{array}{l} v_+ \longmapsto v_+\otimes v_-+v_-\otimes v_+, \\ v_- \longmapsto v_-\otimes v_- \end{array} \right. \\[2mm]
\Delta_2:\quad & W \longrightarrow W\otimes_{\mathbb{Q}}V, & \quad & w \longmapsto w\otimes v_-, \\[2mm]
\theta:\quad & W\otimes_{\mathbb{Q}}W \longrightarrow W\otimes_{\mathbb{Q}}W, & \quad & w\otimes w \longmapsto 0.
\end{aligned}
\]

For each vertex $J$ of the cube, let $M_J(D)$ be the tensor product consisting of one copy of $W$ for each arc and one copy of $V$ for each circle in the state $D(J)$. In other words, each arc and each circle of $D(J)$ contributes a tensor factor of $W$ or $V$ to $M_J(D)$, respectively. Thus, if $D(J)$ consists of $t$ arcs and $r$ circles, then
\[
M_J(D)=
\underbrace{W\otimes_{\mathbb{Q}}\cdots\otimes_{\mathbb{Q}} W}_{t}
\otimes_{\mathbb{Q}}
\underbrace{V\otimes_{\mathbb{Q}}\cdots\otimes_{\mathbb{Q}} V}_{r}.
\]

\subsection{The Khovanov complex}

This cube of $\mathbb{Q}$-vector spaces gives rise to a chain complex. For each integer $i$, the $i$-th chain group is defined as the direct sum of the $\mathbb{Q}$-vector spaces associated with the vertices
\[
\overline{C}^i(D)=\bigoplus_{\substack{J\subseteq I \\ |J|=i}} M_J(D)\{-i\}.
\]
We use $M\{k\}$ to denote the vector space $M$ equipped with a grading shift by $k$. More precisely, if an element $x \in M$ has degree $j$ in $M$, then its degree in $M\{k\}$ is $j-k$. Let $\overline{C}^{i,j}(D)$ be the degree $j$ component of $\overline{C}^i(D)$. Each nonzero saddle has internal degree $-1$, while the state shift $\{-i\}$ raises degree by $i$; hence the differential has bidegree $(1,0)$.

To define the differentials, fix an ordering of the crossings in $I$. Each subset $J \subseteq I$ can then be encoded by a binary sequence
\begin{equation*}
(x_1, x_2,\ldots,x_{|I|})\in\{0,1\}^{|I|},
\end{equation*}
where $x_i=1$ if the $i$-th crossing belongs to $J$, and $x_i=0$ otherwise.

We label each directed edge by an element of $\{0,1,\ast\}^{|I|}$ with exactly one entry equal to $\ast$ and all remaining entries belonging to $\{0,1\}$.
View each labeled edge as a directed arrow whose source is obtained by replacing the symbol $\ast$ by $0$, and whose target is obtained by replacing $\ast$ by $1$. For an edge $\xi \in \{0,1,\ast\}^{|I|}$, let $d_\xi$ denote the corresponding edge map. The map $d_\xi$ acts as the identity on all tensor factors that are unaffected by the corresponding change of resolution. On the factors associated with the components involved in this change, it is given by the appropriate linear map among $m_1$, $m_2$, $\Delta_1$, $\Delta_2$, and $\theta$.

Having defined the edge maps $d_\xi$, we define the differential on a homogeneous element $x \in M_J(D)$ by
\[
d^i(x)
=
\sum_{\substack{\xi}}
(-1)^{\operatorname{sign}(\xi)}d_\xi(x),
\]
where the sum is taken over all directed edges emanating from the vertex $J$. Here, $\operatorname{sign}(\xi)$ is determined by the number of $1$'s appearing before the unique $\ast$ in $\xi$.

\begin{figure}[htbp]
	\centering
	\includegraphics[width=0.25\textwidth]{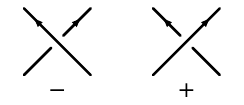}
	\caption{A negative crossing (left) and a positive crossing (right).}
	\label{fig:crossings}
\end{figure}

Taking the orientation of $D$ into account, we obtain the chain complex $C(D)$ from $\overline{C}(D)$ by applying the appropriate grading shifts. The sign for each crossing of $D$ is illustrated in Figure \ref{fig:crossings}.

Let $n_-(D)$ and $n_+(D)$ denote the number of negative and positive crossings in $D$, respectively. The complex $C(D)$ is defined as 
\begin{equation*}
C(D)=\overline{C}(D)[n_-(D)]\{2n_-(D)-n_+(D)\}
\end{equation*}
where the differential remains unchanged. Square brackets indicate a shift in the homological grading. More explicitly,
\[
C^{i,j}(D)=\overline{C}^{\,i+n_-(D),\,j+2n_-(D)-n_+(D)}(D).
\]
An element in $C^{i,j}(D)$ is said to have homological grading $i$ and $q$-grading $j$. If $v \in M_J(D) \subseteq C^{i,*}(D)$ with $q$-grading $j$, then
\begin{equation*}
j=\deg(v)+i+n_+(D)-n_-(D).
\end{equation*}

The statewise comparison in Proposition~\ref{Comparison for every matching} below identifies the signed edge maps on each boundary matching with those of a reduced arc-algebra complex. In particular, $d^2=0$ on $C(D)$. To obtain a tangle invariant, it remains to show that its homology is independent, up to isomorphism, of the choice of diagram $D$. For this purpose, we use Khovanov's arc-algebra construction~\cite{khovanov2002functor}.

Fix $2n$ ordered points on the boundary of a disk, with the orientation signs prescribed by the tangle. Let $\mathscr{T}_n$ denote the set of crossingless matchings of these points, considered up to isotopy relative to the boundary. The matchings themselves are unoriented. We also fix a bending convention $b$ that arranges the boundary points along the boundary line used in Khovanov's construction.

\begin{proposition}\label{Comparison for every matching}
	Let \(D\) be an oriented \(n\)-tangle diagram. For	\(a\in\mathscr T_n\), let \(C_a(D)\) be the space of states with
	boundary matching \(a\). Then
	\[
	C(D)=\bigoplus_{a\in\mathscr T_n}C_a(D)
	\]
	as complexes. Let \(1_b\), \(b\in\mathscr T_n\), be the standard idempotents of the
	nonnegatively graded arc algebra \(H^n\), and set \(P_b=H^n1_b\).
	Following \cite[Section~5]{khovanov2005colored}, define the simple right
	\(H^n\)-module \(\mathbb{Q}_a=\mathbb{Q}u\) by
	\[
	\deg u=0,\qquad u1_b=\delta_{ab}u,\qquad uH^n_{>0}=0.
	\]
	Then there is a chain isomorphism
	\[	C_a(D)\{-n\}\cong
		\mathsf R\!\left(\mathbb{Q}_a\otimes_{H^n}
		\mathcal F_{\mathrm{bnd}}^{\mathrm{Kh}}(bD)\right),\qquad
		\mathcal F_{\mathrm{bnd}}^{\mathrm{Kh}}=\mathcal F^{\mathrm{Kh}}\{n\}_{\mathrm{Kh}}.
	\]
	Here \(\mathsf R\) reverses quantum degree only. The notation \(\{r\}_{\mathrm{Kh}}\) raises the original degree, while the manuscript's
	\(\{r\}\) lowers degree. In particular, \(C(D)\) and each \(C_a(D)\)
	have invariant chain homotopy type.
\end{proposition}

\begin{proof}
	Observe that the differential preserves the boundary matching. Indeed, saddles involving only circles or an arc and a circle do not change the boundary matching. A saddle reconnecting two distinct arcs may change the matching,	but the corresponding map is $\theta=0$. Consequently, the differential preserves the decomposition by boundary matchings, and each $C_a(D)$ is a direct summand subcomplex.
	
	We next compare the state spaces of the two complexes. For a state with matching \(b\) and \(c\) circles, the boundary-normalized arc-algebra complex gives
	\(P_b\otimes A_{\mathrm{Kh}}^{\otimes c}\), where
	\[
	A_{\mathrm{Kh}}=\mathbb{Q}[X]/(X^2),\qquad
	\deg_{\mathrm{Kh}}1=-1,\quad \deg_{\mathrm{Kh}}X=1.
	\]
	Since
	\[
	\mathbb{Q}_a\otimes_{H^n}P_b
	\cong \mathbb{Q}_a1_b
	\cong
	\begin{cases}
	\mathbb{Q}, & a=b,\\
	0, & a\ne b,
	\end{cases}
	\]
	only the states with boundary matching $a$ survive	after tensoring with $\mathbb{Q}_a$.	Reverse the quantum grading. Then $1$ and $X$ are identified with $v_+$ and $v_-$, respectively. The shift $\{-n\}$ comes from the $n$ arc generators
	 $w$, each contributing degree $-1$. Thus the surviving state spaces agree with those of $C_a(D)\{-n\}$.
	
	For a surviving state $s$, denote by $z_{\mathrm{Kh}}$ and $z$ the corresponding
	circle vectors in the two constructions. Their homological and quantum gradings satisfy
	\[
	\begin{aligned}
	i&=|s|-n_-(D),\\
	j&=-\bigl(
	\deg_{\mathrm{Kh}}z_{\mathrm{Kh}}
	-|s|+2n_-(D)-n_+(D)
	\bigr)\\
	&=\deg(z)+|s|+n_+(D)-2n_-(D).
	\end{aligned}
	\]
	These are precisely the gradings of $C_a(D)\{-n\}$.
	
	It remains to compare the differentials. Saddles that involve only circles induce the maps $m_1$ and $\Delta_1$. For a saddle involving an arc and a circle, let $x_\ell$ be the degree-two dot on the corresponding arc in $1_aH^n1_a$. Before tensoring with $\mathbb{Q}_a$, the relevant maps are given by
	\[
	\begin{aligned}
	p\otimes1&\longmapsto p,\\
	p\otimes X&\longmapsto px_\ell,\\
	p&\longmapsto p\otimes X+px_\ell\otimes1.
	\end{aligned}
	\]
	  When tensoring with $\mathbb{Q}_a$,  we find that $x_\ell$ acts trivially. Thus these maps become
	\[
	u\otimes v_+\longmapsto u,\qquad
	u\otimes v_-\longmapsto0,\qquad
	u\longmapsto u\otimes v_-.
	\]
	Identifying \(u\) with the shifted arc generator \(w^{\otimes n}\), we can obtain \(m_2\)
	and \(\Delta_2\) on the affected arc. Note that for a saddle that reconnects two distinct arcs,
	the boundary matching changes. Consequently, at least one of the corresponding state spaces vanishes after the selection. Hence the induced map is zero, which is consistent with $\theta=0$. If the crossing order is the same, then the unnormalized cube signs agree. However, Khovanov's homological shift $[r]$ multiplies the differential by $(-1)^r$ \cite[Sections~2.6 and~3.4]{khovanov2002functor}, while our shift leaves it unchanged. Set $\nu=n_-(D)$. Let $d_{\mathrm{Kh}}$ be the differential on the right-hand side of the claimed isomorphism, and let $\iota_i$ be the statewise identification in normalized homological degree $i$. Then we have
	\[
	d_{\mathrm{Kh}}\iota_i=(-1)^\nu\iota_{i+1}d.
	\]
	Define $\Phi_i=(-1)^{\nu i}\iota_i$. It is clear that 
	\[
	d_{\mathrm{Kh}}\Phi_i
	=(-1)^{\nu(i+1)}\iota_{i+1}d
	=\Phi_{i+1}d,
	\]
Thus the map $\Phi$ is the required chain isomorphism.
	
	Finally, we consider invariance under Reidemeister moves. By~\cite[Theorem~2]{khovanov2002functor}, Khovanov's tangle complexes are related by
	chain homotopy equivalences. Applying $\mathbb{Q}_a\otimes_{H^n}-$ to these equivalences preserves the chain homotopy identities
	$fg-\mathrm{id}=dh+hd$ and the corresponding identities for $gf$. This follows from the additivity of the tensor functor and does not require $\mathbb{Q}_a$ to be flat. The grading shift and reversal of quantum grading also preserve these
	identities. It follows that the selected complexes are invariant up to chain homotopy
	equivalence. Taking their direct sum gives the invariance of the full complex.
	Independence of the crossing order follows from the usual vertex-sign chain isomorphism.
\end{proof}

By Proposition~\ref{Comparison for every matching}, the isomorphism class of the homology of $(C(D),d)$ depends only on the tangle represented by $D$. This allows us to introduce the following definition.

\begin{definition}
	Let $T$ be a tangle and let $D$ be any diagram of $T$. The complex $(C(D),d)$ is called the \textit{Khovanov complex} of $D$, and its homology, up to isomorphism, is called the \textit{Khovanov homology} of $T$, denoted by $H(T)$. We also write $H(D)=H(C(D))$ for the homology of a chosen diagram.
\end{definition}

\subsection{Graded Euler characteristic}

Let $W=\bigoplus_{m\in\mathbb{Z}}W_m$ be a graded vector space over $\mathbb{Q}$. Its \emph{graded (or quantum) dimension} is defined as
\[
\operatorname{qdim}(W)=\sum_{m\in\mathbb{Z}}\dim_\mathbb{Q}(W_m)q^m,
\]
where $\dim_\mathbb{Q}(W_m)$ denotes the dimension of the homogeneous component $W_m$ over $\mathbb{Q}$. Let $C^*=\{C^k\}_{k\in\mathbb{Z}}$ be a chain complex of graded $\mathbb{Q}$-vector spaces. The \emph{graded Euler characteristic} of $C^*$, denoted by $\chi(C^*)$, is defined as
\[
\chi(C^*)=\sum_{k\in\mathbb{Z}}(-1)^k\operatorname{qdim}(C^k).
\]

We first establish a basic property of graded Euler characteristics that will be used in the proof of Theorem \ref{Main}.

\begin{lemma}\label{short exact}
Let $C^{\ast}$ be a bounded complex of finitely generated graded vector spaces over $\mathbb{Q}$. Then 
\[
\chi(C)=\chi(H(C)),
\]
where $H(C)$ is the homology of $C$.
\end{lemma}

\begin{proof}
For each $i$, let
\[
Z^i=\ker(d^i:C^i\to C^{i+1}),
\qquad
B^i=\operatorname{im}(d^{i-1}:C^{i-1}\to C^i).
\]
Since the differential preserves the grading, we have the following two short exact sequences of graded vector spaces
\[
0\longrightarrow Z^i\longrightarrow C^i
\longrightarrow B^{i+1}\longrightarrow 0
\]
and
\[
0\longrightarrow B^i\longrightarrow Z^i
\longrightarrow H^i(C)\longrightarrow 0.
\]
By the additivity of graded dimension, it follows that
\[
\operatorname{qdim}(C^i)
=\operatorname{qdim}(B^i)
+\operatorname{qdim}(H^i(C))
+\operatorname{qdim}(B^{i+1}).
\]
Multiplying the above equation by $(-1)^i$ and summing over all $i$, we get
\[
\sum_i(-1)^i\operatorname{qdim}(C^i)
=\sum_i(-1)^i\operatorname{qdim}(H^i(C))
+\sum_i(-1)^i\big(\operatorname{qdim}(B^i)+\operatorname{qdim}(B^{i+1})\big).
\]
The terms which involves $B^i$ cancel because $C$ is bounded. Hence
$$
\chi(C)
=\sum_i(-1)^i\operatorname{qdim}(C^i)\\
=\sum_i(-1)^i\operatorname{qdim}(H^i(C))\\
=\chi(H(C)).
$$
\end{proof}

This lemma is the graded analogue of the classical fact that the Euler characteristic of a bounded complex equals that of its homology \cite{weibel1994introduction}.

\subsection{Temperley--Lieb theory}\label{sec:tl-theory}

We recall the algebraic and diagrammatic conventions needed in Section~\ref{sec:colored-application}. For further details, see \cite{kauffmanlins1994,blanchet1995topological,khovanov2005colored}. Let
\[
R=\mathbb{Z}[q,q^{-1}],\qquad F=\mathbb{Q}(q),\qquad \delta=q+q^{-1},
\]
and define the quantum integers by
\[
[r]_q=\frac{q^r-q^{-r}}{q-q^{-1}}
=q^{r-1}+q^{r-3}+\cdots+q^{-(r-1)}\quad(r\geq1),
\qquad [0]_q=0.
\]
Then $[1]_q=1$, $[2]_q=\delta$, and $[3]_q=q^2+1+q^{-2}$. Diagram evaluations are carried out over $R$. When rational coefficients are needed for the projectors, we pass to $F$.

For $m\geq 1$, the \emph{Temperley--Lieb algebra} $\operatorname{TL}_m(R,\delta)$ is the unital associative $R$-algebra generated by $e_1,\ldots,e_{m-1}$, which satisfies the relations
\[
e_i^2=\delta e_i,\qquad
e_ie_{i\pm1}e_i=e_i,\qquad
e_ie_j=e_je_i\quad (|i-j|\geq 2),
\]
whenever the indices are defined. The crossingless planar matchings of the $m$ top and $m$ bottom endpoints of a rectangle form an $R$-basis of $\operatorname{TL}_m(R,\delta)$. Two matchings are identified when they are isotopic relative to the boundary. The identity element $I_m$ is represented by $m$ vertical strands, and each generator $e_i$ is represented by a diagram joining the $i$-th and $(i+1)$-st endpoints by a cap at the top and a cup at the bottom, with all remaining strands vertical. Multiplication is defined by vertical stacking, with each closed circle arising from the composition replaced by the scalar $\delta$. For example, stacking $e_i$ with itself produces a single closed circle together with $e_i$, giving $e_i^2=\delta e_i$.

Over $F$, the \emph{Jones--Wenzl projector} $f_m$ is the unique idempotent in
$\operatorname{TL}_m(F,\delta)$ whose coefficient of $I_m$ is $1$ and which satisfies
\[
e_if_m=f_me_i=0
\qquad (1\leq i<m).
\]
In particular, $f_m$ can be written as
$
f_m=I_m+\sum_{a\neq I_m}\lambda_a a,
$
where the sum ranges over the remaining diagram basis elements. The existence and uniqueness of $f_m$ are well known
\cite{kauffmanlins1994,blanchet1995topological}.
For example,
\[
f_1=I_1,\qquad
f_2=I_2-\delta^{-1}e_1.
\]
Using $e_1^2=\delta e_1$, we obtain
\[
f_2e_1
=e_1-\delta^{-1}e_1^2=0.
\]
A further property of the Jones--Wenzl projector, which will be used later, is
given in Lemma~\ref{lem:app-through-ideal}.

Let $V_m$ be the irreducible $U_q(\mathfrak{sl}_2)$-module of highest weight $m$,
with $\dim_F V_m=m+1$. Under the representation-theoretic interpretation of the Temperley--Lieb algebra, the Jones--Wenzl projector $f_m$ projects $V_1^{\otimes m}$
onto its irreducible summand $V_m$. Consequently, an $m$-strand cable carrying $f_m$ is a single strand colored by $N=m+1$, as in Section~\ref{sec:colored-application}.

The  \emph{unnormalised closure trace} is the $F$-linear map
\[
\operatorname{tr}_m:\operatorname{TL}_m(F,\delta)\longrightarrow F
\]
defined as follows. Join each top endpoint of a diagram to the corresponding bottom endpoint outside the rectangle, and assign a factor of $\delta$ to each resulting circle. If  $a$ is a diagram with closure consisting of $\ell(a)$ circles, then
\[
\operatorname{tr}_m(a)=\delta^{\ell(a)}.
\]
In particular,
$
\operatorname{tr}_m(I_m)=\delta^m.
$
Moreover, the closure trace has the property that
\[
\operatorname{tr}_m(xy)=\operatorname{tr}_m(yx)
\]
for all $x,y\in\operatorname{TL}_m(F,\delta)$. 

The closure trace coincides with the representation-theoretic \emph{quantum trace}. If $\mathsf{k}$ is the Cartan operator, then in a weight basis $v_0,\ldots,v_m$ of $V_m$,
\[
\mathsf{k}v_j=q^{m-2j}v_j\qquad(0\leq j\leq m).
\]
We use the pivotal convention in which the loop value is $\delta$. Then the quantum trace of an equivariant endomorphism $A$ of $V_m$ is 
\[
\operatorname{tr}_q(A)=\operatorname{Tr}_{V_m}(\mathsf{k}A).
\]
Since $f_m$ projects onto $V_m$, its closure is the quantum dimension of that summand. Indeed,
\[
\operatorname{tr}_m(f_m)
=\operatorname{tr}_q(\operatorname{id}_{V_m})
=\sum_{j=0}^m q^{m-2j}=[m+1]_q.
\]
More generally, closing a scalar endomorphism multiplies its scalar by the quantum dimension, so
\[
\operatorname{tr}_q(\lambda\operatorname{id}_{V_m})
=\lambda[m+1]_q.
\]
This is the trace factor used in the reduced normalization in Proposition~\ref{prop:app-colored-polynomial}. All projector calculations here are over $F$; specialization at roots of unity is addressed in Remark~\ref{rem:app-root}.

\section{The Jones polynomial of tangles}\label{section3}
In \cite{khovanov2002functor}, Khovanov defined the Kauffman bracket for oriented \((m,n)\)-tangles. In this section, we consider a different class of tangles, namely general \(n\)-tangles, whose definition is not intended as a generalization of that of \((m,n)\)-tangles. We adapt the Kauffman bracket construction to this setting, modifying the signs appearing in the exponents in the definition for oriented \((m,n)\)-tangles to ensure consistency with the Kauffman bracket for links introduced in \cite{khovanov2000categorification}. Based on this modified construction, we define the Jones polynomial for \(n\)-tangles and establish several of its properties. Furthermore, Khovanov proved in \cite{khovanov2000categorification} that the unnormalised Jones polynomial of a link coincides with the graded Euler characteristic of its Khovanov homology. We establish an analogous result for \(n\)-tangles, showing that the scalar invariant \(\mathcal{J}(T)\) associated with an \(n\)-tangle \(T\) is recovered, up to a factor of \(q^n\), as the graded Euler characteristic of the corresponding homology.

\subsection{The Kauffman bracket for tangles}

For an  $n$-tangle diagram $D$, we denote by $\mathcal{R}(D)$ the diagram obtained from $D$ by removing all circles. Let $D$ be a tangle diagram. The \emph{bracket} $\langle D \rangle$ is defined as the following weighted sum
\begin{equation*}
\langle D\rangle=\sum_{D_s}(-q)^{\#(D_s)}(q^{-1}+q)^i \mathcal{R}(D_s),
\end{equation*}
where $D_s$ ranges over the set of all states obtained by resolving $D$, $i$ is the number of circles in $D_s$, and $\#(D_s)$ records the number of $1$-resolutions appearing in the state $D_s$.

\begin{proposition}\label{Kauffman skein}
If $D_1$, $D_2$, and $D_3$ are three $n$-tangle diagrams differing only locally, as shown in Figure \ref{fig:skein_crossings}, then
\begin{equation*}
\langle D_1\rangle=\langle D_2\rangle-q\langle D_3\rangle.
\end{equation*}
\begin{figure}[htbp]
	\centering
	\includegraphics[width=0.5\textwidth]{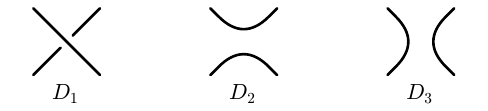}
	\caption{The crossing \(D_1\), its horizontal \(0\)-resolution \(D_2\), and its vertical \(1\)-resolution \(D_3\).}
	\label{fig:skein_crossings}	
\end{figure}	
\end{proposition}	
\begin{proof}
The states of $D_1$ can be partitioned into two classes according to the resolution at the indicated crossing. The states obtained by taking the $0$-resolution at the indicated crossing are in one-to-one correspondence with the states of $D_2$. Since the number of $1$-resolutions is unchanged under this correspondence, their total contribution is $\langle D_2\rangle$.

On the other hand, the states obtained by taking the $1$-resolution at the indicated crossing are in one-to-one correspondence with the states of $D_3$. For each such state, the number of $1$-resolutions increases by one in $D_1$, and hence each state contributes an additional factor of $-q$. Therefore, the total contribution from these states is $-q\langle D_3\rangle$.

Consequently, $\langle D_1\rangle=\langle D_2\rangle-q\langle D_3\rangle$.
\end{proof}

Suppose that $D$ is an oriented diagram of an $n$-tangle $T$. The associated \emph{Kauffman bracket} is defined as
\begin{equation*}
K(D)=(-1)^{n_-(D)}q^{n_+(D)-2n_-(D)}\langle D\rangle.
\end{equation*}	
Under the fixed bending convention $b$ and the resulting identification of boundary matchings, our definition satisfies $K(D;q)=K_{\mathrm{Kh}}(bD;q^{-1})$, where $K_{\mathrm{Kh}}$ denotes the bracket in~\cite{khovanov2002functor}. Thus~\cite[Proposition~11]{khovanov2002functor} implies that $K(D)$ is independent of the choice of diagram $D$ representing the oriented tangle $T$. Hence, $K(D)$ defines a module-valued invariant of $T$, which we denote by $K(T)$.

The Kauffman bracket of a tangle diagram $D$ takes values in the free $\mathbb{Z}[q,q^{-1}]$-module generated by the isotopy classes of circle-free flat $n$-tangles.

As in Section~2.3, let $\mathscr{T}_n$ be the finite set of noncrossing matchings of the fixed $2n$ ordered boundary points of the disk, up to isotopy relative to these points. These matchings represent the flat $n$-tangle diagrams without circles.

We denote the ring $\mathbb{Z}[q,q^{-1}]$ by $R$. For any oriented tangle $T$, the invariant $K(T)$ naturally lies in the free $R$-module $R[\mathscr{T}_n]$. Let $\varphi:R[\mathscr{T}_n]\rightarrow R$ denote the augmentation map, that is, the $R$-module homomorphism given by $\varphi(D)=1$ for every basis element $D\in\mathscr{T}_n$. 
Let $\mathcal{J}(T)$ be the polynomial defined by
\begin{equation*}
\mathcal{J}(T)=\varphi(K(T)).
\end{equation*}
Thus, $\mathcal{J}(T)$ is the scalar invariant obtained from $K(T)$ by  applying $\varphi$. In particular, $\mathcal{J}(\emptyset)=1$. 

\begin{definition}
The \emph{(normalised) Jones polynomial} $\widetilde{\mathcal{J}}(T)$ of an oriented tangle $T$ is defined as
\begin{equation*}
\widetilde{\mathcal{J}}(T)=\frac{\mathcal{J}(T)}{q^{-1}+q}.
\end{equation*}
The polynomial $\mathcal{J}(T)$ is called the \emph{unnormalised Jones polynomial}.
\end{definition}

\begin{remark}
\begin{enumerate}[label=($\roman*$)]
\item For a nonempty link \(L\) with \(\ell\) components, let \(V_L(t)\) denote the classical Jones polynomial with \(V_{\bigcirc}=1\) and skein relation
	\(t^{-1}V_{L_+}-tV_{L_-}=(t^{1/2}-t^{-1/2})V_{L_0}\)
	\cite[Remark~2.12]{murakami2010volume}. Substituting \(t^{1/2}=-q\) gives the skein relation of Theorem~\ref{Jones poly}. Since \(V_L(t)\in t^{(\ell-1)/2}\mathbb Z[t^{\pm1}]\), we obtain
	\[
	\widetilde{\mathcal J}(L;q)
	=\left.V_L(t)\right|_{t^{1/2}=-q}
	=(-1)^{\ell-1}\left.V_L(t)\right|_{t^{1/2}=q}.
	\]
	In particular, \(\widetilde{\mathcal J}(K;q)=V_K(q^2)\) for a knot \(K\).
\item For a tangle \(T\), \(\mathcal{J}(T)\) is, in general, not divisible by
\(q^{-1}+q\). Consequently, \(\widetilde{\mathcal{J}}(T)\), which takes values in
\(\mathbb{Q}(q)\), is in general not a Laurent polynomial.
\end{enumerate}
\end{remark}

\subsection{Properties of the Jones polynomial}

A natural question is whether the tangle polynomial $\mathcal{J}(T)$  can distinguish some tangles whose chosen closures have the same Jones
polynomial. We demonstrate this by comparing a flat tangle and a crossed tangle whose closures are isotopic, where the closure is obtained by connecting the upper endpoints to each other and the lower endpoints to each other.

Consider the trivial $2$-tangle $T_1$ consisting of two parallel vertical arcs, and let $T_2$ be the $2$-tangle with a single positive crossing. Denote the vertical and horizontal flat matchings by $D_v$ and $D_h$, respectively. By the definition of the Kauffman bracket, we have
$
K(T_1)=D_v,
$
and hence
$$
\mathcal{J}(T_1)=\varphi(D_v)=1.
$$
\begin{figure}[H]
	\centering
	\includegraphics[width=0.35\textwidth]{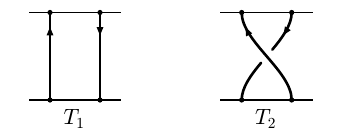}
	\caption{The flat tangle \(T_1\) and the tangle \(T_2\) with one positive crossing.}
	\label{fig:two_tangles}
\end{figure}
For $T_2$, the Kauffman skein relation gives $\langle T_2 \rangle = D_h - q D_v$. Since $n_+(T_2)=1$ and $n_-(T_2)=0$, we have 
\[
K(T_2) = q\langle T_2 \rangle = q D_h - q^{2} D_v.
\]
Applying the augmentation map $\varphi$, we obtain 
\[
\mathcal{J}(T_2) = \varphi(q D_h -q^{2} D_v) = q - q^{2}.
\]
Therefore, we have $\mathcal{J}(T_1) \neq \mathcal{J}(T_2)$. 

However, the numerator closures of $T_1$ and $T_2$ are both isotopic to the unknot and therefore have the same Jones polynomial. Thus, while the classical Jones polynomial cannot distinguish their closures, the tangle invariant $\mathcal{J}$ distinguishes $T_1$ from $T_2$. 

Moreover, the Jones polynomial of tangles also behaves like that of links. 

\begin{theorem} \label{Jones poly}
	The unnormalised Jones polynomial has the following properties:
	\begin{enumerate}[label=($\roman*$)]
		\item If $T$ is a flat tangle without circles, then 
		\begin{equation*}
		\mathcal{J}(T)=1.
		\end{equation*}
		\item There is a skein relation
		\begin{equation*}
		q^{-2}\mathcal{J}(T_1)-q^2\mathcal{J}(T_2)=(q^{-1}-q)\mathcal{J}(T_3),
		\end{equation*}
		whenever $T_1$, $T_2$, and $T_3$ are three oriented tangles that differ only in a local region according to the diagrams in Figure \ref{fig:skein_L1_L2_L3}.
		\begin{figure}[htbp]
			\centering
	        \includegraphics[width=0.5\textwidth]{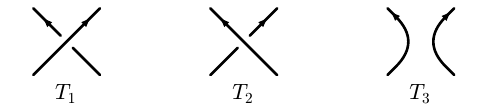}
			\caption{The positive crossing \(T_1\), the negative crossing \(T_2\), and their oriented smoothing \(T_3\).}
			\label{fig:skein_L1_L2_L3}
		\end{figure}
	\end{enumerate}
\end{theorem}

\begin{proof}
$(i)$ By the definition of the Kauffman bracket, $K(T)=\langle T\rangle=(q^{-1}+q)^0 \mathcal{R}(T)=\mathcal{R}(T)=T$. Then $\mathcal{J}(T)=\varphi(K(T))=1$. 

$(ii)$ To establish this skein relation, we need to examine the brackets of the three diagrams $T_1$, $T_2$, and $T_3$. In addition to $T_3$, we introduce another resolution $T_4$, obtained by taking the horizontal resolution of the crossing, while $T_3$ is given by the vertical resolution. Unlike $T_3$, the diagram $T_4$ does not inherit a consistent orientation from the original tangle diagram. Indeed, the orientations of the strands near the resolution would be incompatible, forcing the two orientations to point toward each other. Consequently,  the resolution $T_4$  does not inherit an orientation
compatible with the prescribed exterior orientation, and hence no Jones polynomial can be assigned to it. Nevertheless, the bracket is defined for unoriented diagrams, and we can use $\langle T_4\rangle$ in the following argument.

Note that by Proposition \ref{Kauffman skein},
\begin{align*}
\langle T_1\rangle=\langle T_3\rangle-q\langle T_4\rangle,\\
\langle T_2\rangle=\langle T_4\rangle-q\langle T_3\rangle.
\end{align*}
To cancel the term involving $T_4$, we multiply the first equation by $q^{-1}$ and add the resulting equation to the second one. Thus we obtain
$$q^{-1}\langle T_1\rangle+\langle T_2\rangle=(q^{-1}-q)\langle T_3\rangle.$$
We now apply the definition of $K$ to each of the three diagrams. Since 
\[
\begin{aligned}
n_-(T_1)&=n_-(T_3),\quad &n_-(T_2)&=n_-(T_3)+1,\\
n_+(T_1)&=n_+(T_3)+1,\quad &n_+(T_2)&=n_+(T_3),
\end{aligned}
\]
the above relations yield
\begin{equation*}
q^{-2}K(T_1)-q^2K(T_2)=(q^{-1}-q)K(T_3).
\end{equation*}
Applying \(\varphi\) to both sides of the above equality, we obtain the desired skein relation. This completes the proof.
\end{proof}

\subsection{Categorification of the Jones polynomial for tangles}\label{section4}
Having defined the Jones polynomial for \(n\)-tangles and established its basic properties, we now investigate its relationship with the Khovanov homology introduced in Section 2.3. The following theorem shows that the unnormalised Jones polynomial of an \(n\)-tangle coincides, up to a factor of \(q^n\), with the graded Euler characteristic of its Khovanov homology.

\begin{theorem}\label{Main}
	Let $T$ be an oriented $n$-tangle and let $D$ be any diagram representing $T$.
	Then
	\[
	\mathcal{J}(T)=q^{n}
\chi(H(D)).
	\]
\end{theorem}

\begin{proof}
	Note that for any  finite-dimensional graded $\mathbb{Q}$-vector space $M$,
	the grading shift satisfies
	\[
	\chi(M\{n\})=q^{-n}\chi(M).
	\]
	By Lemma \ref{short exact}, the Euler characteristic of the complex $C(D)$
	can be expressed in terms of its homology:
	\[
	\chi(C(D)) = \chi(H(D)).
	\]
	Therefore, it suffices to prove that
	\begin{equation*}
	\mathcal{J}(T)=q^{n}\chi(C(D)).
	\end{equation*}
	
	Order the distinguished crossing of Figure~\ref{fig:skein_crossings} first. Let
	\(M=\overline{C}(D_2)\) and \(N=\overline{C}(D_3)\{-1\}\) be the complexes of
	the two resolutions, and let \(f\colon M\to N\) be the saddle map
	\cite{bar2005khovanov}. In degree \(i\), both \(\overline{C}(D_1)[1]\) and
	the mapping cone \(\Gamma(f)\) have chain group \(M^{i+1}\oplus N^i\).
	In the \(1\)-resolution the first cube coordinate contributes a minus sign
	to each remaining edge. Since our shift \([1]\) preserves the differential,
	the cube and cone differentials are
	\[
	\begin{aligned}
	d_{[1]}(x,y)&=(d_Mx,\,f(x)-d_Ny),\\
	d_{\Gamma}(x,y)&=(-d_Mx,\,f(x)+d_Ny).
	\end{aligned}
	\]
	Thus \(\Phi_i(x,y)=(-1)^i(x,-y)\) satisfies
	\(d_{\Gamma}\Phi_i=\Phi_{i+1}d_{[1]}\) and gives a chain isomorphism
	\(\overline{C}(D_1)[1]\cong\Gamma(f)\). Taking the alternating sum of the
	graded dimensions of \(M^{i+1}\oplus N^i\) gives
	\(\chi(\Gamma(f))=-\chi(M)+\chi(N)\). Hence
	\[
	\begin{aligned}
	\chi(\overline{C}(D_1))
	&=-\chi(\Gamma(f))=\chi(M)-\chi(N)\\
	&=\chi(\overline{C}(D_2))-q\chi(\overline{C}(D_3)).
	\end{aligned}
	\]

	On the other hand, by Proposition \ref{Kauffman skein}, the quantity $\varphi(\langle D\rangle)$ satisfies the skein relation
\begin{equation*}
\varphi(\langle D_1\rangle)
=
\varphi(\langle D_2\rangle)-q\varphi(\langle D_3\rangle).
\end{equation*}
	Comparing these two skein relations, we see that the quantities $\chi(\overline{C}(D))$ and  $\varphi(\langle D\rangle)$  are characterized by the
	same skein relation. Consequently, for any diagram $D$, repeated applications of
	the skein relation reduce both quantities to the same linear combinations of
	their values on crossingless diagrams. Therefore, to compare these two quantities, it is enough to examine their values on the states of $D$. 
	
	Since every such state is a flat tangle, we only need to determine the values
	of these two quantities on flat tangles. 
 Suppose that $T$ is a flat $n$-tangle with $m$ circles. Then
	\[\overline{C}(T)=
	\underbrace{W\otimes_{\mathbb{Q}}\cdots\otimes_\mathbb{Q} W}_{n}
	\otimes_\mathbb{Q}
	\underbrace{V\otimes_{\mathbb{Q}}\cdots\otimes_\mathbb{Q} V}_{m}.
	\]
Hence
	\[
	\chi(\overline{C}(T))
	=
	q^{-n}(q^{-1}+q)^{m},
	\]
while
$$
	\varphi(\langle T\rangle)
	=\varphi((q^{-1}+q)^{m}\mathcal{R}(T))
	=
	(q^{-1}+q)^{m}.
	$$
	It follows that for any tangle diagram $D$,
\begin{equation*}
\chi(\overline{C}(D))
=q^{-n}\varphi(\langle D\rangle).
\end{equation*}
	
Finally, the relation
\[
C(D)=\overline{C}(D)[n_-(D)]\{2n_-(D)-n_+(D)\}
\]
implies that
\begin{equation*}
\chi(C(D))
=(-1)^{n_-(D)}q^{n_+(D)-2n_-(D)}\chi(\overline{C}(D)).
\end{equation*}
Combining the two preceding identities with the definition of $K(D)$ gives the desired equality $\mathcal{J}(T)=q^n\chi(C(D))$.
\end{proof}

In particular, since a  $0$-tangle is precisely a link, Theorem \ref{Main} recovers the classical relation between the graded Euler characteristic of the Khovanov complex of a link and its Jones polynomial.
\begin{corollary}\cite{bar2002khovanov}
	Let $L$ be an oriented link and let $D$ be any diagram of $L$. The graded Euler characteristic of the complex $C(D)$ is the unnormalized Jones polynomial of $L$.
\end{corollary}

\subsection{Computation and examples}

Two examples are presented next to illustrate Theorem \ref{Main}. For each example, both the Khovanov homology and the Jones polynomial are computed, and the relation between them is verified as stated in Theorem \ref{Main}.

\begin{example}
Consider the Khovanov homology and the Jones polynomial of the Hopf-type tangle $D$ shown in Figure~\ref{fig:hopf-tangle}.
\begin{figure}[H]
	\centering
	\includegraphics[width=0.2\textwidth]{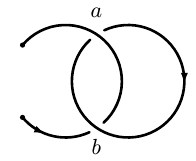}
	\caption{The Hopf-type tangle with negative crossings \(a\) and \(b\).}
	\label{fig:hopf-tangle}
\end{figure}
We first construct a diagram for the complex. Recall that the possible states of the diagram can be organized as the vertices of a $2$-dimensional cube. At each vertex, we place the corresponding diagram obtained from the associated state. This yields Figure \ref{fig:cube-hopf-tangle}.
\begin{figure}[htbp]
	\centering
	\includegraphics[width=0.6\textwidth]{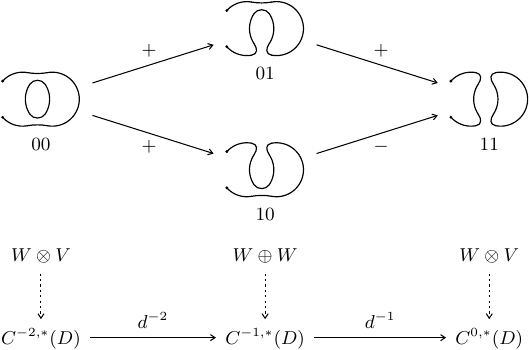}
	\caption{The state cube and Khovanov cochain complex for the Hopf-type tangle. Quantum grading shifts are omitted.}
	\label{fig:cube-hopf-tangle}
\end{figure}

In \(C^{-1,*}(D)\), order the two copies of \(W\) by the states \(10,01\). With this order, the differentials are
\[
\begin{aligned}
d^{-2}(w \otimes v_+)
&=
\bigl(
w,\,
w
\bigr),~ d^{-2}(w \otimes v_-)
=
\bigl(
0,\,
0
\bigr),\\
d^{-1}\bigl((w_1,w_2)\bigr)
&=
-w_1 \otimes v_-+w_2\otimes v_-,
\end{aligned}
\]
where $w_1,w_2\in W$.
To construct the table for this tangle, we need to determine the $q$-gradings. For the diagram under consideration, we have $n_-(D)=2$ and $n_+(D)=0$. The formula for the quantum grading gives the following table.
\[
    \begin{tabular}{l|c|c}
        \toprule
        Height & $-2$ & $0$ \\
        \midrule
        Cycles & $\{w\otimes v_-\}$ & $\{w\otimes v_+,~w\otimes v_-\}$ \\ \hline
        Boundaries & 0 & $\{w\otimes v_-\}$ \\ \hline
        Homology & $\{w\otimes v_-\}$ & $\{w\otimes v_+\}$ \\ \hline
        $q$-degrees & $-6$ & $-2$ \\
        \bottomrule
    \end{tabular}
\]

Then the homology can be summarised as a table where the homological degree is horizontal and the $q$-degree is vertical.
\[
    \begin{tabular}{c|c|c|c}
        \toprule
        $j \backslash i$ & $-2$ & $-1$ & $0$ \\
        \midrule
        $-2$ &  &  & $\mathbb{Q}$ \\ \hline
        $-3$ &  &  &  \\ \hline
        $-4$ &  &  &  \\ \hline
        $-5$ &  &  &  \\ \hline
        $-6$ & $\mathbb{Q}$ &  &  \\
        \bottomrule
    \end{tabular}
\]
It follows that the graded Euler characteristic of the Khovanov homology for this tangle is 
\begin{equation*}
\chi(H(D))=\sum_{i, j\in\mathbb{Z}}(-1)^{i}q^j\operatorname{dim}_{\mathbb{Q}}(H^{i, j}(D))=q^{-2}+q^{-6}.
\end{equation*}

Now we compute the Jones polynomial of this tangle. Let $D^{\prime}$ be the flat 1-tangle.
By the state-sum definition of the bracket, we obtain
\begin{equation*}
\langle D\rangle=(q^{-1}+q)D^{\prime}-2qD^{\prime}+q^2(q^{-1}+q)D^{\prime}=(q^{-1}+q^3)D^{\prime}.
\end{equation*}
Then $K(D)=q^{-4}\langle D\rangle=(q^{-5}+q^{-1})D^{\prime}.$
Thus
\begin{equation*}
\mathcal{J}(D)=q^{-5}+q^{-1}.
\end{equation*}
Comparing the computed Euler characteristic and Jones polynomial gives 
\[\mathcal{J}(D)=q^{1}\chi(H(D)).\]
\end{example}

Having worked out the case of the Hopf-type tangle, which is a 1-tangle with two crossings, we now proceed to a more complicated example, namely, the trefoil-type tangle, which is a 3-tangle consisting of three crossings.
\begin{example}
Consider the trefoil-type tangle illustrated in Figure~\ref{fig:trefoil-tangle}.
\begin{figure}[H]
	\centering
	\includegraphics[width=0.2\textwidth]{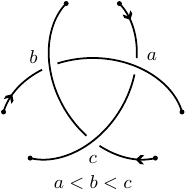}
	\caption{The trefoil-type tangle.}
	\label{fig:trefoil-tangle}
\end{figure}
As in the previous example, we obtain the following complex.
\begin{figure}[htbp]
	\centering
	\includegraphics[width=0.9\textwidth]{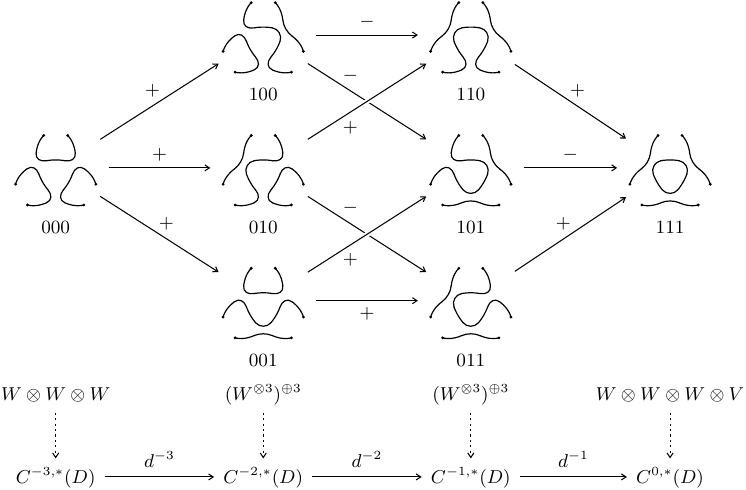}
	\caption{The state cube and Khovanov cochain complex for the trefoil-type tangle. Quantum grading shifts are omitted.}
	\label{fig:cube-trefoil-tangle}
\end{figure}
The differentials in the complex can be written explicitly as
\[
\begin{aligned}
d^{-3}&(w\otimes w\otimes w)
=(0,0,0),\\
d^{-2}&(w_1\otimes w_2\otimes w_3, w_4\otimes w_5\otimes w_6, w_7\otimes w_8\otimes w_9)=(0,0,0),\\
d^{-1}&(w_1\otimes w_2\otimes w_3, w_4\otimes w_5\otimes w_6, w_7\otimes w_8\otimes w_9)=w_1\otimes w_2\otimes w_3\otimes v_-\\
&-w_4\otimes w_5\otimes w_6\otimes v_-+w_7\otimes w_8\otimes w_9\otimes v_-.
\end{aligned}
\]
Here, $w_i\in W$ for $1\leq i\leq 9$. For convenience, we set
$u_i=(w\otimes w\otimes w)_i,~ i=1,2,3.$
Note that the values of $n_-(D)$ and $n_+(D)$ for this diagram $D$ are
$n_-(D)=3, ~n_+(D)=0.$
Taking the quantum gradings into account, we obtain
\[
\begin{tabular}{l|c|c|c|c}
    \toprule
    Height & $-3$ & $-2$ & $-1$ & $0$ \\
    \midrule
    Cycles & $\{w^{\otimes 3}\}$ & $\{u_1, u_2, u_3\}$ & $\left\{\begin{array}{c} u_1+u_2, \\ u_1-u_3 \end{array}\right\}$ & $\left\{\begin{array}{c} w^{\otimes 3}\otimes v_+, \\ w^{\otimes 3}\otimes v_- \end{array}\right\}$ \\
    \hline
    Boundaries & 0& 0& 0& $\{w^{\otimes 3}\otimes v_-\}$ \\        \hline
    Homology & $\{w^{\otimes 3}\}$ & $\{u_1, u_2, u_3\}$ & $\left\{\begin{array}{c} u_1+u_2, \\ u_1-u_3 \end{array}\right\}$ & $\{w^{\otimes 3}\otimes v_+\}$ \\        \hline
    $q$-degrees & $-9$ & $-8$ & $-7$ & $-5$ \\
    \bottomrule
\end{tabular}
\]
Hence, the homology can be summarised as follows.
\[
    \begin{tabular}{c|c|c|c|c}
        \toprule
        $j \backslash i$ & $-3$ & $-2$ & $-1$ & $0$ \\
        \midrule
        $-5$ &  &  &  & $\mathbb{Q}$ \\   \hline
        $-6$ &  &  &  &  \\   \hline
        $-7$ &  &  & $\mathbb{Q}\oplus\mathbb{Q}$ &  \\   \hline
        $-8$ &  & $\mathbb{Q}\oplus\mathbb{Q}\oplus\mathbb{Q}$ &  &  \\   \hline
        $-9$ & $\mathbb{Q}$ &  &  &  \\
        \bottomrule
    \end{tabular}
\]
Then the graded Euler characteristic of the Khovanov homology for this tangle is 
\begin{equation*}
\chi(H(D))=-q^{-9}+3q^{-8}-2q^{-7}+q^{-5}.
\end{equation*}

We now turn to the computation of the Jones polynomial of this tangle.  Applying the state-sum definition of the bracket, we obtain
\begin{equation*}
\varphi(\langle D\rangle)
=1-3q+2q^2-q^4.
\end{equation*}
It follows that 
\begin{equation*}
\mathcal{J}(D)=\varphi(K(D))=-q^{-6}\varphi(\langle D\rangle)=-q^{-6}+3q^{-5}-2q^{-4}+q^{-2}.
\end{equation*}
The computed Euler characteristic and Jones polynomial satisfy
\[
\mathcal{J}(D)=q^3\chi(H(D)).
\]
\end{example}

\begin{figure}[htbp]
	\centering
	\includegraphics[width=0.55\textwidth]{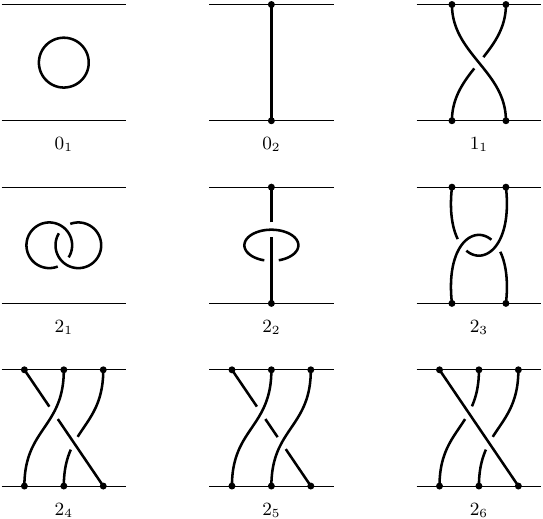}
	\caption{The tangles with at most two crossings listed in Table \ref{computation:at2}.}
	\label{fig:tangle-table-1}
\end{figure}

\begin{figure}[htbp]
	\centering
	\includegraphics[width=0.55\textwidth]{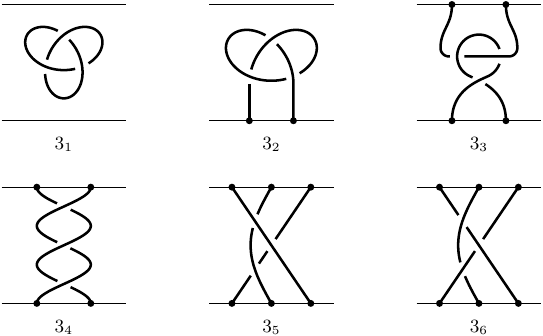}
	\caption{The tangles with three crossings listed in Table \ref{computation:at3}.}
	\label{fig:tangle-table-2}
\end{figure}

Tables~\ref{computation:at2} and~\ref{computation:at3} list the Jones and Poincar\'e polynomials for tangles with at most three crossings. A prime denotes the mirror of the corresponding diagram in Figures~\ref{fig:tangle-table-1} and~\ref{fig:tangle-table-2}. The Poincar\'e polynomials are from \cite{shen2025computing} and give \(\mathcal J(T;q)=q^nP_T(-1,q)\) by Theorem~\ref{Main}.
The last four rows of Table~\ref{computation:at3} cover diagrams consisting of four arcs with a connected projection and three crossings. These tangles are simple in the sense of \cite[Section~3.2]{shen2025computing}. With \(k=n_-(D)\), their common formula \cite[Theorem~2.7]{shen2025computing} is
\[
P_T(t,q)=t^{-k}q^{-1-3k}(1+tq)^3,\qquad 0\leq k\leq3.
\]

\begin{table}[htbp]
	\centering
    \caption{Comparisons of Jones polynomials and Poincar\'e polynomials for tangles with at most two crossings}\label{computation:at2}
	\resizebox{\textwidth}{!}{%
		\begin{tabular}{c|c|c|c|c}
			\toprule
			Tangle type & Sign of crossings & $n$-tangle
			& Poincar\'e polynomial & Unnormalised Jones polynomial \\
			\midrule
			$0_1$ & $\emptyset$ & 0
			& $q+q^{-1}$ & $q+q^{-1}$ \\
			\hline
			$0_2$ & $\emptyset$ & 1
			& $q^{-1}$ & $q\cdot q^{-1}$ \\
			\hline
			\multirow{2}{*}{$1_1$}
			& $\{+\}$ & 2
			& $t+q^{-1}$ & $q^2(-1+q^{-1})$ \\
			& $\{-\}$ & 2
			& $q^{-3}+t^{-1}q^{-4}$ & $q^2(q^{-3}-q^{-4})$ \\
			\hline
			\multirow{2}{*}{$2_1$}
			& $\{+,+\}$ & 0
			& $1+q^2+t^2q^4+t^2q^6$
			& $1+q^2+q^4+q^6$ \\
			& $\{-,-\}$ & 0
			& $1+q^{-2}+t^{-2}q^{-4}+t^{-2}q^{-6}$
			& $1+q^{-2}+q^{-4}+q^{-6}$ \\
			\hline
			\multirow{2}{*}{$2'_1$}
			& $\{+,+\}$ & 0
			& $1+q^2+t^2q^4+t^2q^6$
			& $1+q^2+q^4+q^6$ \\
			& $\{-,-\}$ & 0
			& $1+q^{-2}+t^{-2}q^{-4}+t^{-2}q^{-6}$
			& $1+q^{-2}+q^{-4}+q^{-6}$ \\
			\hline
			\multirow{2}{*}{$2_2$}
			& $\{+,+\}$ & 1
			& $1+t^2q^4$ & $q(1+q^4)$ \\
			& $\{-,-\}$ & 1
			& $q^{-2}+t^{-2}q^{-6}$ & $q(q^{-2}+q^{-6})$ \\
			\hline
			\multirow{2}{*}{$2'_3$}
			& $\{+,+\}$ & 2
			& $1+tq+t^2q^3$ & $q^2(1-q+q^3)$ \\
			& $\{-,-\}$ & 2
			& $q^{-3}+t^{-1}q^{-5}+t^{-2}q^{-6}$
			& $q^2(q^{-3}-q^{-5}+q^{-6})$ \\
			\hline
			\multirow{2}{*}{$2_3$}
			& $\{-,-\}$ & 2
			& $q^{-4}+t^{-1}q^{-5}+t^{-2}q^{-7}$ & $q^2(q^{-4}-q^{-5}+q^{-7})$ \\
			& $\{+,+\}$ & 2
			& $q^{-1}+tq+t^2q^{2}$
			& $q^2(q^{-1}-q+q^{2})$ \\
			\hline
			\multirow{3}{*}{$2_4$, $2_5$, $2_6$}
			& $\{+,+\}$ & 3
			& $t^2q+2t+q^{-1}$
			& $q^3(q-2+q^{-1})$ \\
			& $\{+,-\}$ & 3
			& $tq^{-2}+2q^{-3}+t^{-1}q^{-4}$
			& $q^3(-q^{-2}+2q^{-3}-q^{-4})$ \\
			& $\{-,-\}$ & 3
			& $q^{-5}+2t^{-1}q^{-6}+t^{-2}q^{-7}$
			& $q^3(q^{-5}-2q^{-6}+q^{-7})$ \\
			\bottomrule
		\end{tabular}%
	}
\end{table}

\begin{table}[htbp]
	\centering
    \caption{Comparisons of Jones polynomials and Poincar\'e polynomials for tangles with three crossings}\label{computation:at3}
	\begin{adjustbox}{max width=\textwidth}
		\begin{tabular}{c|c|c|c|c}
			\toprule
			Tangle & Crossing signs & $n$-tangle
			& Poincar\'e polynomial & Unnormalised Jones polynomial \\
			\midrule
			$3_1$ & $\{-,-,-\}$ & 0
			& $t^{-3}q^{-9}+t^{-2}q^{-5}+q^{-3}+q^{-1}$
			& $-q^{-9}+q^{-5}+q^{-3}+q^{-1}$ \\
			\hline
			$3'_1$ & $\{+,+,+\}$ & 0
			& $t^3q^9+t^2q^5+q^3+q$
			& $-q^9+q^5+q^3+q$ \\
			\hline
			$3_2$ & $\{-,-,-\}$ & 1
			& $t^{-3}q^{-9}+t^{-2}q^{-7}+q^{-3}$
			& $q(-q^{-9}+q^{-7}+q^{-3})$ \\
			\hline
			$3'_2$ & $\{+,+,+\}$ & 1
			& $t^3q^7+t^2q^5+q$
			& $q(-q^7+q^5+q)$ \\
			\hline
			\multirow{2}{*}{$3_3$}
			& $\{-,-,-\}$ & 2
			& $t^{-3}q^{-9}+t^{-2}q^{-8}+t^{-2}q^{-7}+t^{-1}q^{-6}+q^{-4}$
			& $q^2(-q^{-9}+q^{-8}+q^{-7}-q^{-6}+q^{-4})$ \\
			& $\{+,+,-\}$ & 2
			& $t^2q^2+t+q^{-1}+q^{-2}+t^{-1}q^{-3}$
			& $q^2(q^2-1+q^{-1}+q^{-2}-q^{-3})$ \\
			\hline
			\multirow{2}{*}{$3'_3$}
			& $\{+,+,+\}$ & 2
			& $t^3q^5+t^2q^4+t^2q^3+tq^2+1$
			& $q^2(-q^5+q^4+q^3-q^2+1)$ \\
			& $\{-,-,+\}$ & 2
			& $t^{-2}q^{-6}+t^{-1}q^{-4}+q^{-3}+q^{-2}+tq^{-1}$
			& $q^2(q^{-6}-q^{-4}+q^{-3}+q^{-2}-q^{-1})$ \\
			\hline
			\multirow{2}{*}{$3_4$}
			& $\{+,+,+\}$ & 2
			& $t^3q^6+t^2q^4+tq^2+q$
			& $q^2(-q^6+q^4-q^2+q)$ \\
			& $\{-,-,-\}$ & 2
			& $t^{-3}q^{-8}+t^{-2}q^{-7}+t^{-1}q^{-5}+q^{-3}$
			& $q^2(-q^{-8}+q^{-7}-q^{-5}+q^{-3})$ \\
			\hline
			\multirow{2}{*}{$3'_4$}
			& $\{-,-,-\}$ & 2
			& $t^{-3}q^{-10}+t^{-2}q^{-8}+t^{-1}q^{-6}+q^{-5}$
			& $q^2(-q^{-10}+q^{-8}-q^{-6}+q^{-5})$ \\
			& $\{+,+,+\}$ & 2
			& $t^3q^4+t^2q^3+tq+q^{-1}$
			& $q^2(-q^4+q^3-q+q^{-1})$ \\
			\hline
			\multirow{2}{*}{$3'_5$}
			& $\{+,+,+\}$ & 3
			& $2t^2q^2+2tq+1$
			& $q^3(2q^2-2q+1)$ \\
			& $\{+,-,-\}$ & 3
			& $t^{-2}q^{-6}+2t^{-1}q^{-5}+2q^{-4}$
			& $q^3(q^{-6}-2q^{-5}+2q^{-4})$ \\
			\hline
			\multirow{2}{*}{$3_5$}
			& $\{-,-,-\}$ & 3
			& $2t^{-2}q^{-8}+2t^{-1}q^{-7}+q^{-6}$
			& $q^3(2q^{-8}-2q^{-7}+q^{-6})$ \\
			& $\{-,+,+\}$ & 3
			& $t^2+2tq^{-1}+2q^{-2}$
			& $q^3(1-2q^{-1}+2q^{-2})$ \\
			\hline
			\multirow{2}{*}{$3_6$}
			& $\{-,-,-\}$ & 3
			& $t^{-3}q^{-9}+3t^{-2}q^{-8}+2t^{-1}q^{-7}+q^{-5}$
			& $q^3(-q^{-9}+3q^{-8}-2q^{-7}+q^{-5})$ \\
			& $\{+,+,-\}$ & 3
			& $t^2q+2tq^{-1}+3q^{-2}+t^{-1}q^{-3}$
			& $q^3(q-2q^{-1}+3q^{-2}-q^{-3})$ \\
			\hline
			\multirow{2}{*}{$3'_6$}
			& $\{+,+,+\}$ & 3
			& $t^3q^3+3t^2q^2+2tq+q^{-1}$
			& $q^3(-q^3+3q^2-2q+q^{-1})$ \\
			& $\{-,-,+\}$ & 3
			& $tq^{-3}+3q^{-4}+2t^{-1}q^{-5}+t^{-2}q^{-7}$
			& $q^3(-q^{-3}+3q^{-4}-2q^{-5}+q^{-7})$ \\
			\hline
			\multirow{4}{*}{Simple $4$-tangles}
			& $\{+,+,+\}$ & 4
			& $q^{-1}+3t+3t^2q+t^3q^2$
			& $q^4(q^{-1}-3+3q-q^2)$ \\
			& $\{+,+,-\}$ & 4
			& $t^{-1}q^{-4}+3q^{-3}+3tq^{-2}+t^2q^{-1}$
			& $q^4(-q^{-4}+3q^{-3}-3q^{-2}+q^{-1})$ \\
			& $\{+,-,-\}$ & 4
			& $t^{-2}q^{-7}+3t^{-1}q^{-6}+3q^{-5}+tq^{-4}$
			& $q^4(q^{-7}-3q^{-6}+3q^{-5}-q^{-4})$ \\
			& $\{-,-,-\}$ & 4
			& $t^{-3}q^{-10}+3t^{-2}q^{-9}+3t^{-1}q^{-8}+q^{-7}$
			& $q^4(-q^{-10}+3q^{-9}-3q^{-8}+q^{-7})$ \\
			\bottomrule
		\end{tabular}
	\end{adjustbox}
\end{table}

\section{Disjoint union, gluing, and mirror symmetry}\label{section5}

Having defined the Jones polynomial via the projection of the module-valued Kauffman bracket, we now examine how it behaves under basic topological operations.
In this section we consider three operations, namely disjoint union, gluing of $1$-tangles along the boundary, and mirror reflection.

\subsection{Disjoint union of tangles}

Unlike the gluing operation, disjoint union does not identify any boundary points, so no closed circles are created. The scalar invariant is therefore multiplicative under disjoint union, with no extra loop factors.

\begin{theorem}\label{thm:tangle-disjoint-union}
	Let $T_1$ and $T_2$ be oriented $n_1$-tangle and $n_2$-tangle, respectively. Then
	\begin{equation*}
	\mathcal{J}(T_1\sqcup T_2)
	=\mathcal{J}(T_1)\mathcal{J}(T_2).
	\end{equation*}
\end{theorem}

\begin{proof}
	Let $D_1$ and $D_2$ be diagrams representing $T_1$ and $T_2$, respectively. Each state $s$ of $D_1\sqcup D_2$ is uniquely determined by a pair of states $(s_1,s_2)$ of $D_1$ and $D_2$. Since
	\[\begin{aligned}
	n_{\pm}(D_1\sqcup D_2)
	&=n_{\pm}(D_1)+n_{\pm}(D_2),\\
	\#(s)&=\#(s_1)+\#(s_2),\\
	i(s)&=i(s_1)+i(s_2),
	\end{aligned}
	\]
	the corresponding state weights satisfy
	\[	w(s)=w(s_1)w(s_2).\]
	Furthermore,
$\mathcal{R}(s)=\mathcal{R}(s_1)\sqcup\mathcal{R}(s_2).$
	
	Consider the $R$-linear map induced by disjoint union,
	\[\mu:
	R[\mathscr{T}_{n_1}]\otimes_R
	R[\mathscr{T}_{n_2}]
	\longrightarrow R[\mathscr{T}_{n_1+n_2}],
	\]
	defined on basis elements by
	\[\mu(D\otimes D')
	=D\sqcup D'.
	\]
	It follows from the state-sum definition that
$$K(T_1\sqcup T_2)
	=\sum_{s_1,s_2}w(s_1)w(s_2)
	\bigl(\mathcal{R}(s_1)\sqcup
	\mathcal{R}(s_2)\bigr)\\
	=\mu\bigl(K(T_1)\otimes K(T_2)\bigr).$$
	
	Since the augmentation sends every basis element to $1$, we have
	$\varphi(D\sqcup D')=\varphi(D)\varphi(D')=1.$
	By linearity, this gives
	\[
	\varphi\bigl(\mu(x\otimes y)\bigr)
	=\varphi(x)\varphi(y)
	\]
	for all $x\in R[\mathscr{T}_{n_1}]$ and $y\in R[\mathscr{T}_{n_2}]$. Consequently,
	$$\mathcal{J}(T_1\sqcup T_2)
	=\varphi\bigl(K(T_1\sqcup T_2)\bigr)\\
	=\varphi(K(T_1))\varphi(K(T_2))\\
	=\mathcal{J}(T_1)\mathcal{J}(T_2).$$
\end{proof}

\subsection{Gluing formula for \texorpdfstring{$1$}{1}-tangles}

In this subsection, we establish a gluing formula for the unnormalised Jones polynomial of two $1$-tangles. Since $R[\mathscr{T}_1]$ is a free $R$-module of rank one, the module-valued invariant $K(T)$ is determined by the scalar invariant $\mathcal{J}(T)$. This allows us to express the Jones polynomial of the resulting link in terms of the invariants of the two tangles.

\begin{theorem}\label{thm:tangle_gluing}
Let $T_1$ and $T_2$ be oriented $1$-tangles with compatible orientations at the boundary
points to be identified. Let $T$ be the oriented link obtained by gluing the upper and lower endpoints of $T_1$ to the lower and upper endpoints of $T_2$, respectively, without
introducing additional crossings. Then
\[
\mathcal{J}(T)
=(q^{-1}+q)\mathcal{J}(T_1)\mathcal{J}(T_2).
\]
\end{theorem}

\begin{proof}
For $n=1$, any circle-free flat tangle consists of a single arc joining the two boundary points. Since all such arcs are isotopic relative to the boundary, $\mathscr{T}_1$ contains only one element, denoted by $D_1$. Hence,
\[
R[\mathscr{T}_1]=R\cdot D_1\cong R.
\]

For any $1$-tangle $T$, the module-valued invariant $K(T) \in R[\mathscr{T}_1]$ can therefore be expressed as a scalar polynomial multiplied by this single basis element. Since the $R$-module homomorphism $\varphi: R[\mathscr{T}_n] \to R$ is defined by $\varphi(D)=1$ for all $D \in \mathscr{T}_n$, we have:
\begin{equation*}
K(T) = \mathcal{J}(T) \cdot D_1.
\end{equation*}
This implies that for $1$-tangles, the module-valued invariant $K(T)$ and the scalar invariant $\mathcal{J}(T)$ carry identical information.

Consider the gluing $T=T_1\cup T_2$. Since the gluing introduces no additional crossings, each state $s$ of $T$ corresponds uniquely to a pair of states $(s_1,s_2)$ of $T_1$ and $T_2$. For $i=1,2$, the state $s_i$ consists of $i(s_i)$ circles and a single arc isotopic to $D_1$. Gluing these two arcs produces one additional circle. Hence,
\[
i(s)=i(s_1)+i(s_2)+1.
\]
After removing all closed components, the resulting flat tangle is the empty tangle $\emptyset\in\mathscr{T}_0$.

Since the gluing does not change the crossings of either tangle, the crossing numbers and the numbers of $1$-resolutions are additive. Together with the preceding identity, we have
\[
w(s)=(q^{-1}+q)w(s_1)w(s_2).
\]
If we sum over all pairs of states $(s_1, s_2)$, then
\begin{equation*}
K(T) = \sum_{s_1, s_2} w(s_1) w(s_2) (q^{-1}+q) \, \emptyset.
\end{equation*}
Note that the choices of states $s_1$ for $T_1$ and $s_2$ for $T_2$ are independent. Then the double sum factors into the product of two sums, which means
\begin{equation*}
K(T) = (q^{-1}+q) \left( \sum_{s_1} w(s_1) \right) \left( \sum_{s_2} w(s_2) \right) \emptyset = (q^{-1}+q) \mathcal{J}(T_1) \mathcal{J}(T_2) \, \emptyset.
\end{equation*}

Finally, the map \(\varphi\) sends the tangle invariant to the scalar invariant of the closed link. Since $\varphi(\emptyset) = 1$, 
\begin{equation*}
\mathcal{J}(T) = \varphi(K(T)) = (q^{-1}+q) \, \mathcal{J}(T_1) \, \mathcal{J}(T_2),
\end{equation*}
which is the desired result.
\end{proof}

\subsection{Reflection and mirror symmetry}

Another fundamental topological operation is taking the mirror image of a tangle, which reverses all crossings. The behavior of the Jones polynomial under this operation is a classic property. We now verify that the Jones polynomial $\mathcal{J}(T)$ satisfies the same mirror relation.

\begin{theorem}\label{thm:tangle-mirror}
Let $T$ be an oriented $n$-tangle, and let $T^*$ denote its mirror image obtained by reversing all crossings. Then
\begin{equation*}
\mathcal{J}(T^*)(q)
=\mathcal{J}(T)(q^{-1}).
\end{equation*}
\end{theorem}

\begin{proof}
Let $D$ be a diagram of $T$, and let $D^*$ be the corresponding mirror diagram.
Reversing all crossings interchanges positive and negative crossings, so
\[n_+(D^*)=n_-(D),\qquad
n_-(D^*)=n_+(D).
\]

There is a bijection between the states of $D$ and those of $D^*$ obtained by
interchanging the two resolutions at each crossing. If $s$ and $s^*$ are corresponding
states, their resolutions have the same underlying flat tangle. In particular,
\[
i(s^*)=i(s).
\]
It follows that the factors $(q^{-1}+q)^{i(s)}$ are unchanged under this correspondence.
Let\[
p=n_+(D),\qquad
m=n_-(D),\qquad
c=p+m.
\] If a state \(s\) has \(r\) \(1\)-resolutions and \(\ell\) circles, then its total weight is
\[
w_D(s)=(-1)^{m+r}q^{p-2m+r}(q^{-1}+q)^\ell.
\]
Since the corresponding state \(s^*\) of the mirror diagram \(D^*\) has \(c-r\) \(1\)-resolutions,
\[
\begin{aligned}
w_{D^*}(s^*)
&=(-1)^{p+c-r}q^{m-2p+c-r}(q^{-1}+q)^\ell\\
&=(-1)^{m+r}q^{-(p-2m+r)}(q^{-1}+q)^\ell\\
&=w_D(s)\big|_{q\mapsto q^{-1}}.
\end{aligned}
\]
So the weight of the corresponding mirror state is $w_D(s)$ with $q$ replaced by $q^{-1}$.
It is obvious that switching a crossing interchanges its $0$- and $1$-resolutions. Then the complementary states $s$ and $s^*$ therefore give the same planar matching, so $\mathcal{R}(s^*)=\mathcal{R}(s)$.
Hence the module-valued invariant of the mirror tangle is
\begin{equation*}
K(T^*) = \sum_{s^*} w_{T^*}(s^*) \mathcal{R}(s^*) = K(T)\big|_{q \mapsto q^{-1}}.
\end{equation*}

Finally, applying the map $\varphi$,  we obtain 
\begin{equation*}
\mathcal{J}(T^*)(q) = \varphi\left( K(T)\big|_{q \mapsto q^{-1}} \right) = \varphi(K(T))(q^{-1}) = \mathcal{J}(T)(q^{-1}).
\end{equation*}
This completes the proof.
\end{proof}

\section{Application to reduced colored invariants and closure volume}\label{sec:colored-application}

Applying the augmentation $\varphi$ to $K(T)$ gives the scalar $\mathcal{J}(T)$. This sums the coefficients
of all boundary matchings. Here we keep the coefficient of one matching in a parallel cable. Throughout this
section, let $T$ be a zero-framed oriented $(1,1)$-tangle with a single arc, fixed endpoints, and no closed
components. Its standard closure is a knot $\widehat{T}$. We use the color $N=m+1\geq2$ as in
Section~\ref{sec:tl-theory}.

\subsection{The colored through coefficient}

On the diagram basis of Section~\ref{sec:tl-theory}, define $\epsilon_m:\operatorname{TL}_m(R,\delta)\to R$
by
\[
\epsilon_m(I_m)=1,\qquad \epsilon_m(a)=0\quad (a\neq I_m)
\]
and extend it $R$-linearly. We call a state $s$ of $D$ \emph{through} if $\mathcal{R}(D_s)=I_m$. This means
that each arc joins a bottom endpoint to a top endpoint. The state may also contain circles.

\begin{lemma}\label{lem:app-through-ideal}
The map $\epsilon_m$ is a unital algebra homomorphism. The Jones--Wenzl projector $f_m$ satisfies
\begin{equation*}
f_mxf_m=\epsilon_m(x)f_m,\qquad x\in\operatorname{TL}_m(\mathbb{Q}(q),\delta).
\end{equation*}
\end{lemma}
\begin{proof}
See \cite[Section~3]{blanchet1995topological}.
\end{proof}

Choose a diagram $D$ of $T$ with blackboard framing zero. For a crossing-free closing arc, this means that
$D$ has writhe zero. We can arrange this by adding curls. Let $D_m$ be the blackboard $m$-parallel cable and
orient its strands coherently. Define
\begin{equation*}
\mathcal{J}_N(T;q)=\epsilon_m(K(D_m)),\qquad m=N-1.
\end{equation*}
Thus $\mathcal{J}_N(T;q)$ is the coefficient of the identity matching in the bracket of the cable.

\begin{proposition}\label{prop:app-colored-polynomial}
Let $J_N(L;t)$ denote the reduced colored Jones polynomial of a zero-framed knot. It takes the value $1$ on
the unknot. With $t=q^2$, we have
\begin{equation*}
\mathcal{J}_N(T;q)=J_N(\widehat{T};q^2).
\end{equation*}
In particular, \[
\mathcal{J}_2(T;q)=\mathcal{J}(T),
\qquad
\mathcal{J}_N(U;q)=1,\]
where $U$ denotes the trivial arc.
\end{proposition}

\begin{proof}
	Let $m=N-1$. Consider the zero-framed $m$-cable $D_m$ of $T$. We write
	\[ K(D_m)=c_{I_m}(q)I_m+
	\sum_{a\neq I_m}c_a(q)a,
	\]
	where \[
	c_{I_m}(q)
	=\epsilon_m(K(D_m))
	=\mathcal{J}_N(T;q).
	\]
	It is important to note that only through states contribute to this coefficient. The state sum gives
	\[
	c_{I_m}(q)=
	\sum_{\substack{s\,:\,
			\mathcal{R}((D_m)_s)=I_m}}
	(-1)^{n_-+|s|}
	q^{n_+-2n_-+|s|}
	(q+q^{-1})^{c(s)},
	\]
	where $n_\pm=n_\pm(D_m)$, $|s|$ is the number of $1$-resolutions, and $c(s)$
	is the number of circles in the state. In particular,	$c_{I_m}(q)\in\mathbb{Z}[q,q^{-1}].	$
	
	By Lemma~\ref{lem:app-through-ideal},
	\[
	f_mK(D_m)f_m=c_{I_m}(q)f_m.
	\]
	To obtain the reduced colored Jones polynomial, we take the closure trace and divide by
	$[N]_q$~\cite{khovanov2005colored}. For our zero-framed diagram, we have
	\[w(D_m)=m^2w(D)=0.
	\]
	It follows that the framing correction is trivial. The components of the closed cable have pairwise linking number zero.
	Thus the alternating orientations used in~\cite{khovanov2005colored} give the same normalized bracket.
	
	Since $\operatorname{tr}_m(f_m)=[N]_q$, we obtain
	\[
	\begin{aligned}
	J_N(\widehat{T};q^2)
	&=\frac{\operatorname{tr}_m
		(f_mK(D_m)f_m)}{[N]_q}\\
	&=\frac{c_{I_m}(q)
		\operatorname{tr}_m(f_m)}{[N]_q}\\
	&=c_{I_m}(q)
	=\mathcal{J}_N(T;q).
	\end{aligned}
	\]
	
	For $N=2$, we have $m=1$, $D_1=D$, and $\epsilon_1=\varphi$. Therefore,
	\[
	\mathcal{J}_2(T;q)
	=\epsilon_1(K(D))
	=\varphi(K(D))
	=\mathcal{J}(T).
	\]
	Finally, the $m$-cable of the trivial arc is $I_m$, and hence
	\[
	\mathcal{J}_N(U;q)
	=\epsilon_m(I_m)=1.
	\]
\end{proof}

\begin{remark}\label{rem:app-root}
	We first work over $F$ with generic $q$, where the Jones--Wenzl projector is defined. The state sum gives
	$\mathcal{J}_N(T;q)$ in $\mathbb{Z}[q,q^{-1}]$. In the volume calculation below, we evaluate this Laurent
	polynomial at $q=e^{\pi i/N}$. At this value $[N]_q$ vanishes, so we take the trace quotient before
	specializing.
\end{remark}

\subsection{The reduced colored complex and comparison}

Let $B$ be an oriented tangle diagram with $m$ bottom and $m$ top endpoints. A circle saddle or an
arc--circle saddle preserves the boundary matching. A saddle that reconnects two distinct arcs has edge map
$\theta=0$. Thus the differential of $C(B)$ preserves each matching. The through states form a direct
summand subcomplex of $C(B)$, denoted by $C_I(B)$. Define
\begin{equation*}
C_{\mathrm{th}}(B)=C_I(B)\{-m\},
\qquad
C_N(T;D)=C_{\mathrm{th}}(D_m).
\end{equation*}
By the convention of Section~\ref{sec:preliminaries}, the shift $\{-m\}$ raises quantum degree by $m$. This
cancels the total degree $-m$ of the generators of the $m$ arc factors $W$. Omitting these generators leaves
a copy of $V^{\otimes c(s)}$ for each through state $s$ with $c(s)$ circles. A homogeneous element $z$ has
bigrading
\begin{equation*}
i=|s|-n_-(B),\qquad
j=\deg(z)+|s|+n_+(B)-2n_-(B).
\end{equation*}

We use the differential of the cube. The maps on circle factors are $m_1$ and $\Delta_1$. If we identify the
distinguished arc generators with $1$, the arc--circle maps become
\begin{equation*}
1\otimes v_+\longmapsto 1,\qquad
1\otimes v_-\longmapsto 0,\qquad
1\longmapsto 1\otimes v_-.
\end{equation*}
The cube signs are those of $C(B)$. All edge maps between through and non-through states vanish. The
following comparison gives the invariance of $C_{\mathrm{th}}(B)$.

\begin{theorem}\label{thm:app-comparison}
	Let $T$ be an oriented $1$-tangle as above and let $N=m+1$. The chain homotopy type of $C_N(T;D)$ is
	independent of the zero-framed diagram $D$ and the order of its crossings. In our gradings, its homology
	$H_N(T)$ is isomorphic to Khovanov's finite reduced $m$-colored homology of $\widehat{T}$. Its graded Euler
	characteristic satisfies
	\begin{equation*}
	\chi(H_N(T))
	=\mathcal{J}_N(T;q)
	=J_N(\widehat{T};q^2).
	\end{equation*}
	For $N=2$, we have $C_2(T;D)=C(D)\{-1\}$. The identity then gives the single-arc case of
	Theorem~\ref{Main}.
\end{theorem}

\begin{proof}
	Bending the identity tangle gives the crossingless matching $e=b(I_m)$. Let $\mathbb{Q}_e$ be the
	corresponding simple right $H^m$-module. By Proposition~\ref{Comparison for every matching}, there is a chain
	isomorphism
	\[
	C_{\mathrm{th}}(B)\xrightarrow{\;\Phi_B\;}
	\mathsf R\!\left(
	\mathbb{Q}_e\otimes_{H^m}
	\mathcal F^{\mathrm{Kh}}_{\mathrm{bnd}}(bB)
	\right).
	\]
	In homological degree $i$, the map is given by $(\Phi_B)_i=(-1)^{n_-(B)i}\iota_i$. Here $\iota_i$ is the statewise
	identification in that proposition. The complex $\mathcal F^{\mathrm{Kh}}_{\mathrm{bnd}}(bB)$ is Khovanov's arc-algebra complex~\cite{khovanov2002functor}
	with the boundary grading shift $\{m\}_{\mathrm{Kh}}$. Tensoring with $\mathbb{Q}_e$ over $H^m$ gives
	Khovanov's simple-module reduction~\cite[Section~5]{khovanov2005colored}.
	
	Reidemeister moves give chain homotopy equivalences by~\cite[Theorem~2]{khovanov2002functor}. Tensoring
	with $\mathbb{Q}_e$ preserves these equivalences. Indeed, applying the tensor functor to
	\[fg-\mathrm{id}=dh+hd
	\]
	gives the same homotopy identity on the reduced complexes. Shifting the gradings and reversing quantum
	degree also preserve this identity. Diagrams of the same zero-framed tangle give isotopic framed cables, so
	the reduced complexes have the same chain homotopy type. A change of crossing order gives the usual chain
	isomorphism by a change of signs at the vertices.
	
	Khovanov uses alternating orientations on the cable~\cite[Section~5]{khovanov2005colored}. Let $S$ index
	the components of the closed cable whose orientations we reverse. Only crossings between a reversed
	component and a component whose orientation is fixed change sign. The components have pairwise linking
	number zero by the framing assumption. Hence the change in writhe is
	\[
	w'-w
	=-4\sum_{\substack{r\in S\\t\notin S}}
	\operatorname{lk}(L_r,L_t)=0.
	\]
	The total number of crossings stays the same, so $n'_\pm=n_\pm$. Both orientations give the same cube and
	the same grading shifts. We thus obtain Khovanov's finite reduced colored complex, with the unknot in
	bidegree $(0,0)$.
	
	Write $n_\pm=n_\pm(B)$. Summing the graded dimensions of the through states with their homological signs
	gives
	\[
	\begin{aligned}
	\chi(C_{\mathrm{th}}(B))
	&=\sum_{\mathcal R(B_s)=I_m}
	(-1)^{|s|-n_-}
	q^{|s|+n_+-2n_-}
	(q+q^{-1})^{c(s)}\\
	&=\epsilon_m(K(B)).
	\end{aligned}
	\]
	Taking $B=D_m$ and using Proposition~\ref{prop:app-colored-polynomial},
	we obtain
	$$\chi(H_N(T))
	=\chi(C_N(T;D))\\
	=\mathcal{J}_N(T;q)\\
	=J_N(\widehat{T};q^2).
	$$
	
	For $N=2$, we have $m=1$ and $D_1=D$. Every state of a $1$-tangle has the unique
	boundary matching, so
	\[
	C_2(T;D)=C(D)\{-1\}.
	\]
	This completes the proof.
\end{proof}

\begin{remark}
	A resolution is called backtracking if it has an arc with both endpoints on the same side. Non-backtracking
	states also occur in sutured Khovanov homology~\cite{grigsbywehrli2010,grigsbyni2014}. Here we use the
	signed even complex over $\mathbb{Q}$. The local $V,W$ construction realizes Khovanov's finite reduced
	theory defined by simple modules. The comparison applies only to this theory.
\end{remark}

\subsection{Stacking and colored gluing}

Let $T$ and $S$ be two tangles with compatible boundary orientations and framings. Denote their vertical composition by $T\circ S$. Then the closure of $T\circ S$ is isotopic to the connected sum $\widehat{T}\#\widehat{S}$.

\begin{theorem}\label{thm:app-stacking}
	Let $T$ and $S$ be zero-framed single-arc tangles with compatible boundary orientations. Then
	\begin{align*}
	C_N(T\circ S)
	&\simeq C_N(T)\otimes_{\mathbb{Q}}C_N(S),\\
	H_N(T\circ S)
	&\cong H_N(T)\otimes_{\mathbb{Q}}H_N(S).
	\end{align*}
	In the stacked diagram $B_1\circ B_2$, order the crossings of $B_1$ before those of $B_2$. On a through
	state $(s_1,s_2)$, the chain isomorphism identifies the circle factors and multiplies by the sign
	\[
	(-1)^{n_-(B_1)|s_2|}.
	\]
\end{theorem}

\begin{proof}
	Choose states $s_1$ and $s_2$ of $B_1$ and $B_2$. Denote their circle-free matchings by $a_1$ and $a_2$,
	respectively. By Lemma~\ref{lem:app-through-ideal},
	\[
	\epsilon_m(a_1a_2)
	=\epsilon_m(a_1)\epsilon_m(a_2).
	\]
	The composite state is therefore through exactly when $a_1=a_2=I_m$. Thus a through state of $B_1\circ B_2$
	is a pair of through states $(s_1,s_2)$. Stacking two identity matchings produces no circles, so
	\[
	c(s_1,s_2)=c(s_1)+c(s_2).
	\]
	The corresponding state spaces are identified by
	\[
	V^{\otimes(c(s_1)+c(s_2))}
	\cong
	V^{\otimes c(s_1)}
	\otimes_{\mathbb Q}
	V^{\otimes c(s_2)}.
	\]
	Moreover,
	\[
	|s|=|s_1|+|s_2|,
	\qquad
	n_\pm(B_1\circ B_2)
	=n_\pm(B_1)+n_\pm(B_2),
	\]
	so the normalized homological and quantum gradings are additive.
	
	Let $d_1$ and $d_2$ be the differentials of the two factors, and let $\nu=n_-(B_1)$. The identification of circle factors preserves the saddle maps. With the chosen crossing order, the cube differential is
	\[
	d_{\mathrm{stack}}
	=d_1\otimes 1+(-1)^{|s_1|}1\otimes d_2,
	\]
	and the tensor-product differential is
	\[
	d_{\otimes}
	=d_1\otimes 1+(-1)^{|s_1|-\nu}1\otimes d_2.
	\]
	Define $\Phi$ on the summand corresponding to $(s_1,s_2)$ by
	\[
	\Phi(z_1\otimes z_2)
	=(-1)^{\nu|s_2|}z_1\otimes z_2.
	\]
	An edge of $B_1$ leaves the exponent $\nu|s_2|$ unchanged. For an edge of $B_2$, the signs agree because
	\[
	(-1)^{|s_1|+\nu(|s_2|+1)}
	=
	(-1)^{\nu|s_2|+|s_1|-\nu}.
	\]
	So
	$
	\Phi d_{\mathrm{stack}}
	=d_{\otimes}\Phi
	$
	and $\Phi$ is a chain isomorphism
	\[
	C_{\mathrm{th}}(B_1\circ B_2)
	\cong
	C_{\mathrm{th}}(B_1)
	\otimes_{\mathbb Q}
	C_{\mathrm{th}}(B_2).
	\]
	Cabling commutes with stacking for zero-framed diagrams. Applying Theorem~\ref{thm:app-comparison} gives
	the first assertion and the K\"unneth theorem over $\mathbb Q$ gives
	\[
	H_N^{i,j}(T\circ S)
	\cong
	\bigoplus_{\substack{i_1+i_2=i\\
			j_1+j_2=j}}
	H_N^{i_1,j_1}(T)
	\otimes_{\mathbb Q}
	H_N^{i_2,j_2}(S),
	\]
	which proves the second assertion.
\end{proof}

\begin{corollary}\label{cor:app-gluing}
	Let $T$ and $S$ be zero-framed single-arc tangles with compatible boundary orientations. The Poincaré polynomial of $H_N$ and its graded Euler characteristic are multiplicative under stacking. In particular,
	\begin{equation*}
	\mathcal{J}_N(T\circ S;q)
	=\mathcal{J}_N(T;q)\mathcal{J}_N(S;q).
	\end{equation*}
	Taking closures with compatible framings gives
	\begin{align*}
	J_N(\widehat{T}\#\widehat{S};q^2)
	&=\mathcal{J}_N(T;q)\mathcal{J}_N(S;q),\\
	J_N^{\mathrm{unred}}(\widehat{T}\#\widehat{S};q^2)
	&=[N]_q\mathcal{J}_N(T;q)\mathcal{J}_N(S;q).
	\end{align*}
	For $N=2$, the second identity reduces to the gluing formula in Theorem~\ref{thm:tangle_gluing}.
\end{corollary}

\begin{proof}
	Let
	$P_N(T;t,q)
	=\sum_{i,j}t^iq^j
	\dim_{\mathbb Q}H_N^{i,j}(T).$
	By Theorem~\ref{thm:app-stacking},
	\[
	P_N(T\circ S;t,q)
	=P_N(T;t,q)P_N(S;t,q).
	\]
	Set $t=-1$. By Theorem~\ref{thm:app-comparison},
	\[
	\mathcal{J}_N(T\circ S;q)
	=\mathcal{J}_N(T;q)\mathcal{J}_N(S;q).
	\]
	As a framed knot, the closure of $T\circ S$ is isotopic to
	$\widehat{T}\#\widehat{S}$.
	The first colored Jones identity follows from Proposition~\ref{prop:app-colored-polynomial}.
	The second follows from the normalization
	\[
	J_N^{\mathrm{unred}}(L;q^2)
	=[N]_qJ_N(L;q^2).
	\]
	Finally, for $N=2$, we have $[2]_q=q+q^{-1}$ and
	$\mathcal{J}_2=\mathcal{J}$, which gives the formula in Theorem~\ref{thm:tangle_gluing}.
\end{proof}

\begin{remark}
Corollary~\ref{cor:app-gluing} uses the Jones--Wenzl projector on the closed cable. Closing the undecorated
cable instead gives the ordinary Jones polynomial. Theorem~\ref{thm:app-stacking} applies to vertical stacking only, and its formula does not hold as stated for arbitrary colored tangle cobordisms.
\end{remark}

\subsection{The figure-eight tangle and stacking families}

For a knot $L$, let $\operatorname{Vol}(L)$ be the sum of the volumes of the hyperbolic pieces in the JSJ
decomposition of its exterior. If \(L\) is hyperbolic, \(\operatorname{Vol}(L)\) is its usual hyperbolic volume. It vanishes for torus knots and is additive under connected sum. It also equals $v_3$ times the
relative simplicial volume of the exterior, where $v_3$ is the volume of a regular ideal
tetrahedron~\cite{murakamimurakami2001,murakami2010volume}.

For a tangle $T$, define
\[
\operatorname{Vol}_{\mathrm{cl}}(T)
=\operatorname{Vol}(\widehat{T}),
\qquad
q_N=e^{\pi i/N}.
\]
Here $\operatorname{Vol}_{\mathrm{cl}}(T)$ measures the volume of the closure $\widehat{T}$.
By Theorem~\ref{thm:app-comparison}, the identity
\begin{equation*}
\lim_{N\to\infty}\frac{2\pi}{N}
\log\left|
\left.\chi(H_N(T))\right|_{q=q_N}
\right|
=\operatorname{Vol}_{\mathrm{cl}}(T)
\end{equation*}
is equivalent to the reduced colored Jones volume conjecture for $\widehat{T}$. We use the convention
$\log 0=-\infty$.

\begin{proposition}\label{prop:app-figure-eight}
	Let $T_E$ be a zero-framed single-arc tangle whose closure is the figure-eight knot $E=4_1$.
	Then
	\[
	\lim_{N\to\infty}\frac{2\pi}{N}
	\log\left|
	\left.\chi(H_N(T_E))\right|_{q=q_N}
	\right|
	=2v_3.
	\]
\end{proposition}

\begin{proof}
	By Theorem~\ref{thm:app-comparison},
	\[
	\left.\chi(H_N(T_E))\right|_{q=q_N}
	=J_N(E;e^{2\pi i/N}).
	\]
	For the figure-eight knot, the logarithmic limit of this colored Jones evaluation is
	$\operatorname{Vol}(E)=2v_3$~\cite[Section~3.2]{murakami2010volume}.
\end{proof}

\begin{corollary}\label{cor:app-volume-family}
	Let $T$ be a finite stack of zero-framed single-arc tangles. Suppose that $r$ factors come from cutting
	figure-eight knots. Each remaining factor comes from cutting a torus knot or the unknot.  Mirrors of these factors are also allowed.
	Then
	\begin{equation*}
	\lim_{N\to\infty}\frac{2\pi}{N}
	\log\left|
	\left.\chi(H_N(T))\right|_{q=q_N}
	\right|
	=2rv_3=\operatorname{Vol}_{\mathrm{cl}}(T).
	\end{equation*}
\end{corollary}

\begin{proof}
	Let $T=T_1\circ\cdots\circ T_k$. By Corollary~\ref{cor:app-gluing} and
	Theorem~\ref{thm:app-comparison},
	\[
	\left.\chi(H_N(T))\right|_{q=q_N}
	=\prod_{j=1}^{k}
	J_N(\widehat{T_j};e^{2\pi i/N}).
	\]
	The limit for a figure-eight factor is $2v_3$~\cite[Section~3.2]{murakami2010volume}. For a torus knot it
	is zero~\cite[Corollary in Section~1]{kashaevtirkkonen2000}.
	The reduced colored Jones polynomial of the unknot is $1$. For a mirror image, the substitution $q\mapsto q^{-1}$ leaves the absolute value unchanged at $q=q_N$, since
	$q_N^{-1}=\overline{q_N}$.
	
	Taking logarithms of the product and applying these limits gives
	\[
	\lim_{N\to\infty}\frac{2\pi}{N}
	\log\left|
	\left.\chi(H_N(T))\right|_{q=q_N}
	\right|
	=2rv_3.
	\]
	The closure of $T$ is the connected sum of the knots $\widehat{T_j}$. Volume is additive under connected
	sum, so
	\[
	\operatorname{Vol}_{\mathrm{cl}}(T)
	=\sum_{j=1}^{k}\operatorname{Vol}(\widehat{T_j})
	=2rv_3.
	\]
	This proves the assertion.
\end{proof}

\begin{remark}
	These limits follow from known cases of the volume conjecture for knots. For a composite closure, we sum
	the volumes of the hyperbolic pieces of its exterior.
	
	We also obtain a lower bound on the total dimension of the homology. Since $|q_N|=1$, the triangle
	inequality gives
	\[
	\left|
	\left.\chi(H_N(T))\right|_{q=q_N}
	\right|
	\leq \dim_{\mathbb{Q}}H_N(T).
	\]
	It follows that
	\[
	\liminf_{N\to\infty}
	\frac{2\pi}{N}
	\log\dim_{\mathbb{Q}}H_N(T)
	\geq 2rv_3.
	\]
	Terms in the graded Euler characteristic can cancel, so this argument gives only a lower bound.
\end{remark}

\section*{Acknowledgments}

This work was supported in part by the National Natural Science Foundation of China (NSFC Grant No. 12401080), the Science and Technology Research Program of Chongqing Municipal Education Commission (Grant No. KJQN202601151), the Scientific Research Foundation of Chongqing University of Technology, the National Key Research and Development Program of China (NKPs Grant No. 2024YFA1013201), the Natural Science Foundation of Chongqing (NSFCQ Grant No. CSTB2024NSCQ-LZX0040), and the Special Project of Chongqing Municipal Science and Technology Bureau (Grant No. 2025CCZ015). Li Shen was supported by an NITMB fellowship supported by grants from the Simons Foundation (Grant No. MP-TMPS-00005320).

\section*{Declaration of AI assistance}

In writing this paper, the authors used ChatGPT to polish the language, double-check the proofs, and format the manuscript. All mathematical statements, proofs, and conclusions were checked independently by the authors, who take full responsibility for the paper.

\bibliographystyle{plain}  % unsrtnat
\bibliography{Reference}

\end{document}